# Linearization and Lemma of Newton for Operator functions

M. Stiefenhofer

**Abstract.** We study the action of the nonlinear mapping $G[z]$ between real or complex Banach spaces in the vicinity of a given curve with respect to possible linearization, emerging patterns of level sets, as well as existing solutions of $G[z] = 0$. The results represent local generalizations of the standard implicit or inverse function theorem and of Newton's Lemma, considering the order of approximation needed to obtain solutions of $G[z] = 0$.

The main technical tool is given by Jordan chains with increasing rank, used to obtain an Ansatz, appropriate for transformation of the nonlinear system to its linear part. The family of linear mappings is restricted to the case of an isolated singularity.

Geometrically, the Jordan chains define a generalized cone around the given curve, composed of approximate solutions of order $2k$ with $k$ denoting the maximal rank of Jordan chains needed to ensure $k$-surjectivity of the linear family. Along these lines, the zero set of $G[z]$ in the cone is calculated immediately, agreeing up to the order of $k-1$ with the given approximation. Hence, the results may also be interpreted as a version of Tougeron's implicit function theorem or Hensel's Lemma in Banach spaces, essentially restricted to the arc case of a single variable.

Finally, by considering a left shift of the Jordan chains, the Ansatz can be modified in a systematic way to obtain a sequence of refined versions of linearization theorems and Newton Lemmas in Banach spaces.

*Contents*

## *1. Introduction and Main Results*

Given a smooth mapping $G[\cdot] \in C^r(B, \bar{B})$ with $B, \bar{B}$ real or complex Banach spaces and a smooth curve $z(\cdot) \in C^r(U, B), 0 \in U \subset \mathbb{K} = \mathbb{R}, \mathbb{C}$ satisfying $z(0) = 0$ and $r \geq 2$ chosen sufficiently large.

*Matthias Stiefenhofer, University of Applied Sciences Kempten, 87435 Kempten (Germany)*
*matthias.stiefenhofer@hs-kempten.de*

We study the action of the mapping $G[z]$ from a vicinity of the curve $z(\varepsilon)$ in $B$ to a vicinity of the curve $G[z(\varepsilon)]$ in $\bar{B}$ with respect to possible linearization, emerging patterns of level sets, as well as existing solutions of $G[z] = 0$. The results represent local generalizations of the standard implicit or inverse function theorem and of Newton's Lemma, considering the order of approximation needed to obtain solutions of $G[z] = 0$.

First, Taylor's formula along $z(\varepsilon)$ implies for $b \in B$

$$\begin{aligned} G[\,z(\varepsilon)+b\,] \;&=\; G[\,z(\varepsilon)\,] \;+\; G'[\,z(\varepsilon)\,]\cdot b \;+\; \tfrac{1}{2}G''[\,z(\varepsilon)\,]\cdot b^2 \;+\; r(\varepsilon,b)\cdot b^3 \\ &=: \; G[\,z(\varepsilon)\,] \;+\; L(\varepsilon)\cdot b \;+\; \mathcal{B}(\varepsilon)\cdot b^2 \;+\; r(\varepsilon,b)\cdot b^3 \end{aligned} \tag{1.1}$$

with $G[z(\varepsilon)]$ describing the behaviour of $G[z]$ along the curve $z(\varepsilon)$, $L(\cdot) \in C^{r-1}(U, L[B,\bar{B}])$ defining the directional derivative and $\mathcal{B}(\cdot) \in C^{r-2}(U, L^2[B,\bar{B}])$ giving the directional curvature at $z(\varepsilon)$. Here $L[B,\bar{B}]$ and $L^2[B,\bar{B}]$ denote bounded linear and bilinear operators from $B$ to $\bar{B}$ respectively. The mapping $r(\cdot)$ represents a smooth remainder function.

Now, if the family of linear operators $L(\varepsilon) = G'[0] + \varepsilon \cdot r_1(\varepsilon) = L_0 + \varepsilon \cdot r_1(\varepsilon)$ has a constant term $L_0$, which is surjective, then, using the decomposition

$$B = N_1^c \;\oplus\; N_1 \tag{1.2}$$

with nullspace $N_1 \coloneqq N[L_0]$, we obtain by standard implicit function theorem (assume $N_1^c$ and $N_1$ to be closed) the existence of a local smooth transformation in $B$ of the form

$$b = \underbrace{\psi_1^c(\,\varepsilon, \varphi_1, n_1\,)}_{\to N_1^c} + n_1\,, \quad \varphi_1 \in N_1^c\,, \;\; n_1 \in N_1\,, \tag{1.3}$$

eliminating higher order terms in (1.1) according to

$$G[\,z(\varepsilon) + \psi_1^c(\,\varepsilon, \varphi_1, n_1\,) + n_1\,] = G[\,z(\varepsilon)\,] \;+\; L_0 \cdot \varphi_1 \tag{1.4}$$

and $\psi_1^c(\cdot)$ satisfying at $\varepsilon = 0$ the near identity relation

$$\psi_1^c(\,0, \varphi_1, n_1) + n_1 = (\,\varphi_1 + n_1\,) \;+\; \chi_1^c(n_1, \varphi_1)\cdot \begin{pmatrix} n_1 \\ \varphi_1 \end{pmatrix}^2. \tag{1.5}$$

Upper Index $c$ in (1.3)-(1.5) means image space of the mapping is given by the complement $N_1^c$ of the kernel $N_1$.

Further, if the constant term $G[0]$ of $G[z(\varepsilon)] = G[0] + \varepsilon \cdot r_2(\varepsilon)$ is sufficiently small, e.g. $G[0] = 0$, then $G[z(\varepsilon)]$ in (1.4) may also be absorbed into the linear part $L_0$ by a small shift of $\varphi_1 \in N_1^c$ according to

$$G[\,z(\varepsilon) + \psi_1^c(\,\varepsilon,\; \varphi_1 - L_0^{-1}\cdot G[z(\varepsilon)]\,,\; n_1\,) + n_1\,] = L_0 \cdot \varphi_1 \tag{1.6}$$

and complete linearization is achieved.

From (1.4) we see that for fixed values of $(\varepsilon, \varphi_1)$, the right hand side of the equation turns into the constant value $G[z(\varepsilon)] + L_0 \cdot \varphi_1$ in $\bar{B}$, and thus, the corresponding level set to this value is given by $b = z(\varepsilon) + \psi_1^c(\varepsilon, \varphi_1, n_1) + n_1$, implying each level set near $z(\varepsilon)$ to be parametrized by $n_1 \in N_1$.

In addition, choosing $(\varepsilon, \varphi_1) = (0,0)$ in (1.6), solutions of $G[z] = 0$ near $z(0) = 0$ are given by $b = \psi_1^c(0, -L_0^{-1} \cdot G[0], n_1) + n_1$.

To obtain the partial linearization (1.4), the surjectivity condition $R[L_0] = \bar{B}$ is needed with respect to the family $L(\varepsilon)$, which has to be supplemented by a smallness requirement with respect to $G[z(\varepsilon)]$ to end up with complete linearization (1.6).

Our aim is to generalize this procedure concerning linearization and calculation of level and zero sets. First, the surjectivity condition with respect to the family $L(\varepsilon)$ generalizes as follows.

There exists $k \geq 0$, such that for every $\bar{b}$ in the image space $\bar{B}$, there exists a curve

$$b(\varepsilon) = b_0 + \varepsilon \cdot b_1 + \cdots + \varepsilon^k \cdot b_k \tag{1.7}$$

in the domain $B$ satisfying

$$L(\varepsilon) \cdot b(\varepsilon) = \varepsilon^k \cdot \bar{b} + \varepsilon^{k+1} \cdot r(\varepsilon)\,, \tag{1.8}$$

i.e. the curve $b(\varepsilon)$ represents an approximation of order $k$ with respect to the linear equation $L(\varepsilon) \cdot b = 0$ with leading coefficient given by $\varepsilon^k \cdot \bar{b}$.

In case of $k = 0$, condition (1.8) turns into

$$L(\varepsilon) \cdot b(\varepsilon) = (\, L_0 + \varepsilon \cdot L_1 + \cdots) \cdot b_0 = L_0 \cdot b_0 + \varepsilon \cdot r(\varepsilon) \overset{(1.8)}{=} \bar{b} + \varepsilon \cdot r(\varepsilon) \tag{1.9}$$

and for each $\bar{b} \in \bar{B}$ we are requested to find $b_0$ with $L_0 \cdot b_0 = \bar{b}$. Hence, if $k = 0$, condition (1.8) simplifies to the surjectivity condition of the standard implicit function theorem.

Note that condition (1.8) merely depends on the Taylor polynomial $L_k(\varepsilon) \coloneqq L_0 + \cdots + \varepsilon^k \cdot L_k$ of degree $k$ of $L(\varepsilon)$ and in this sense, the directional derivative operators $L(\varepsilon)$ are supposed to be $k$-determined with respect to surjectivity, as defined by (1.7), (1.8).

Note also that for $\varepsilon \neq 0$, the range of the linear operator $L(\varepsilon)$ satisfies $R[L(\varepsilon)] = \bar{B}$ by (1.8), implying $L(\varepsilon)$ to be a surjective family of linear operators, with the exception of $L(0) = L_0$ allowed to be a singular operator, i.e. at $\varepsilon = 0$ an isolated singularity may exist. We call a family $L(\varepsilon)$ of linear operators satisfying (1.8) a $k$-surjective family. The minimal number $k \geq 0$ satisfying (1.8) is called the order of surjectivity of the family $L(\varepsilon)$.

If $L(0)$ does not turn into a surjective operator for $\varepsilon \neq 0$, then $k$-surjectivity of $L(\varepsilon)$ may be generalized to stabilization at $k$ of the family $L(\varepsilon)$. Here we refer to [2], [10], [18].

For checking condition (1.8) with $k \geq 1$, Jordan chains of increasing rank have to be calculated for $L(\varepsilon)$ until surjectivity of leading coefficients occurs in $\bar{B}$. During this process direct sums of $B$ and $\bar{B}$ are build up according to the following scheme.

$$\begin{array}{ccccccccccc} & & \overbrace{\phantom{N}}^{rank=0} & & \overbrace{\phantom{N}}^{rank=1} & & & & \overbrace{\phantom{N}}^{rank=k} & & \\ B & = & N_1^c & \oplus & N_2^c & \oplus & \cdots & \oplus & N_{k+1}^c & \oplus & N_{k+1} \\ & & \boxed{\updownarrow S_1} & & \boxed{\updownarrow S_2} & & & & \boxed{\updownarrow S_{k+1}} & & \\ \bar{B} & = & R_1 & \oplus & R_2 & \oplus & \cdots & \oplus & R_{k+1} & & \end{array} \tag{1.10}$$

Here, the subspace $N_i^c, i = 1, \cdots, k+1$ contains root elements of Jordan chains of rank $i - 1$ with corresponding leading coefficients in $\bar{B}$ defining the subspace $R_i$. The correlation between ele-ments in $N_i^c$ and leading coefficients in $R_i$ is established by an operator $S_i \in L[B,\bar{B}]$, representing an isomorphism $S_i \in GL[N_i^c, R_i]$ when restricted to $N_i^c$. Finally, $N_{k+1}$ contains root elements that can be continued up to arbitrary high order of approximation with respect to the linear equation $L(\varepsilon)\cdot b = 0$.

Note also that for $k = 0$, the split condition (1.10) turns into decomposition (1.2) of the standard implicit function theorem under consideration of $S_1 = L_0 = G'[0]$.

Now, assuming $k$-surjectivity of the family $L(\varepsilon)$, the linear part $L(\varepsilon)$ in (1.1) allows for direct impact on every element in $\bar{B}$, ultimately permitting complete elimination of higher order terms in (1.1) according to the following linearization theorem.

**Theorem 1:** Assume the family $L(\varepsilon)$ to be $k$-surjective with $k \geq 0$ and all subspaces closed in the direct sums of (1.10).

(i) Then a smooth transformation in $B$ of the form

$$b = p_k(\varepsilon)\cdot\{\underbrace{\Psi^c(\varepsilon,\varphi,n_{k+1})}_{\to N_1^c\oplus\cdots\oplus N_{k+1}^c} + n_{k+1}\} \tag{1.11}$$

$$p_k(\varepsilon) = \underbrace{\phi_0}_{=I_B} + \varepsilon\cdot\phi_1 + \cdots + \varepsilon^k\cdot\phi_k\,, \quad \phi_i \in L[B, N_1^c\oplus\cdots\oplus N_{k+1-i}^c]\,,\ \ i = 1,\ldots,k \tag{1.12}$$

$$\varphi = (\varphi_1,\cdots,\varphi_{k+1}) \in N^c := N_1^c\times\cdots\times N_{k+1}^c\,, \quad n_{k+1}\in N_{k+1} \tag{1.13}$$

exists, eliminating higher order terms in (1.1) according to

$$\begin{aligned} &G[z(\varepsilon) + \varepsilon^k\cdot p_k(\varepsilon)\cdot\{\Psi^c(\varepsilon,\varphi,n_{k+1}) + n_{k+1}\}] \\ &= G[z(\varepsilon)] + \varepsilon^{2k}\cdot(S_1\cdots S_{k+1})\cdot\varphi \end{aligned} \tag{1.14}$$

with $(S_1\cdots S_{k+1}) \in GL[N^c,\bar{B}]$ and triangularity of the operator polynomial

$$P_k(\varepsilon) := S_1 + \cdots + \varepsilon^k\cdot S_{k+1}\,.$$

In addition, the nonlinear mapping $\Psi^c(\cdot)$ satisfies

$$\Psi^c(\varepsilon,\varphi,n_{k+1}) = \varepsilon^k\cdot\underbrace{\psi_1^c(\varepsilon,\varphi,n_{k+1})}_{\to N_1^c} + \cdots + \varepsilon^1\cdot\underbrace{\psi_k^c(\varepsilon,\varphi,n_{k+1})}_{\to N_k^c} + \underbrace{\psi_{k+1}^c(\varepsilon,\varphi,n_{k+1})}_{\to N_{k+1}^c}, \tag{1.15}$$

as well as

$$\begin{array}{ccccc} \psi_1^c(0,\varphi,n_{k+1}) & = & \varphi_1 & + & \chi_1^c(n_{k+1},\varphi_{k+1})\cdot\begin{pmatrix} n_{k+1} \\ \varphi_{k+1}\end{pmatrix}^2 \\ \vdots & & \vdots & & \vdots \\ \psi_k^c(0,\varphi,n_{k+1}) & = & \varphi_k & + & \chi_k^c(n_{k+1},\varphi_{k+1})\cdot\begin{pmatrix} n_{k+1} \\ \varphi_{k+1}\end{pmatrix}^2 \\ \psi_{k+1}^c(0,\varphi,n_{k+1}) & = & \varphi_{k+1} & + & \chi_{k+1}^c(n_{k+1},\varphi_{k+1})\cdot\begin{pmatrix} n_{k+1} \\ \varphi_{k+1}\end{pmatrix}^2 \end{array} \tag{1.16}$$

at $\varepsilon = 0$ and $\psi_1^c(\varepsilon,0,0) = \cdots = \psi_{k+1}^c(\varepsilon,0,0) = 0$ for $(\varphi,n_{k+1}) = 0$.

(ii) For fixed $\varepsilon \in U, \varepsilon \neq 0$, the mapping

$$z = z(\varepsilon) + \varepsilon^k \cdot p_k(\varepsilon) \cdot \{ \Psi^c(\varepsilon, \varphi, n_{k+1}) + n_{k+1} \} \tag{1.17}$$

represents a diffeomorphism between

$$V = \{ (\varphi, n_{k+1}) \in N^c \times N_{k+1} \mid \|(\varphi, n_{k+1})\| \ll 1 \} \tag{1.18}$$

and a local neighbourhood of the curve position $z(\varepsilon)$ in $B$.

Finally, for fixed values of $(\varepsilon, \varphi), \varepsilon \neq 0$, the level set to the value

$$\bar{b} = G[\, z(\varepsilon)\, ] + \varepsilon^{2k} \cdot (S_1 \cdots S_{k+1}) \cdot \varphi \tag{1.19}$$

is given by (1.17), parametrized by $n_{k+1} \in N_{k+1}$.

In case of $k = 0$, formulas (1.11)-(1.19) are again simplifying to the standard case of (1.3)-(1.5) by neglecting $\varphi_1, \cdots, \varphi_k$ and replacing $(\varphi_{k+1}, n_{k+1})$ by $(\varphi_1, n_1)$.

From a more geometrical viewpoint, the mapping (1.17) represents an $\varepsilon$-blow up transformation in $B$, guided by expansion rates $\varepsilon^0, \varepsilon^1, \cdots, \varepsilon^k$, which are ensured by (1.8) and (1.10) to be present within the family $L(\varepsilon)$ of directional derivatives along $z(\varepsilon)$. These expansion rates are calculated via Jordan chains of $L(\varepsilon)$ and in this sense, Theorem 1 mainly takes into consideration the linear part $L(\varepsilon) \cdot b$ of the Taylor expansion (1.1). However, when looking in detail for solutions of the equation $G[z] = 0$, then the first term $G[z(\varepsilon)]$ and the higher order terms $\mathcal{B}(\varepsilon) \cdot b^2 + r(\varepsilon, b) \cdot b^3$ in (1.1) have to be considered too.

In figure 1, the working principle of the $\varepsilon$-parametric family of diffeomorphisms (1.17) is indicated in a schematic way in $\mathbb{R}^n$.

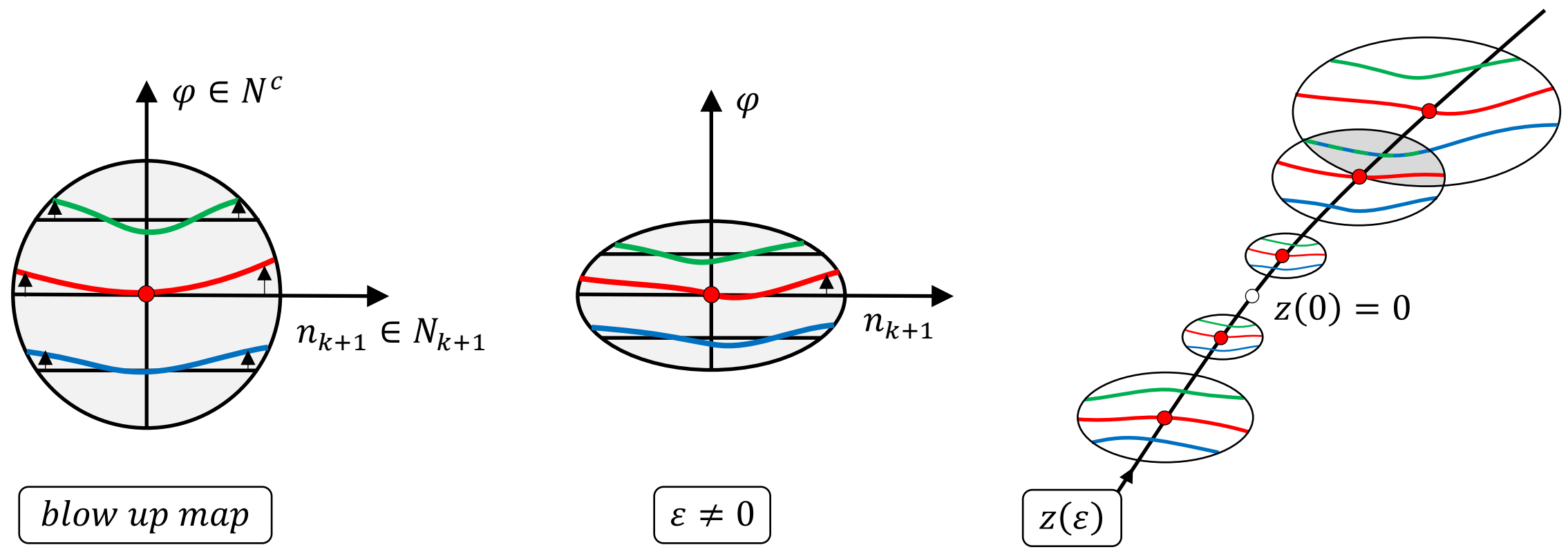

Figure 1: Schematic interpretation of the blow up procedure.

The mapping is constructed in several steps starting on the left for $\varepsilon = 0$ with the sum operator

$$(\varphi, n_{k+1}) \;\rightarrow\; \psi_1^c(0, \varphi, n_{k+1}) \;+\; \cdots \;+\; \psi_k^c(0, \varphi, n_{k+1}) \;+\; \psi_{k+1}^c(0, \varphi, n_{k+1})\,, \tag{1.21}$$

which is bending each horizontal line $\varphi = const$ smoothly into the direction of $N^c$ (indicated by arrows), thereby approximating the identity map towards $(\varphi, n_{k+1}) \rightarrow 0$ according to (1.16). In addition, the red marked image of the $N_{k+1}$ axis is tangentially touching the $N_{k+1}$ axis in the origin. The map in the left diagram represents the blow up map.

Next, when choosing $\varepsilon \neq 0$ according to

$$(\varphi, n_{k+1}) \;\rightarrow\; \psi_1^c(\varepsilon, \varphi, n_{k+1}) + \cdots + \psi_k^c(\varepsilon, \varphi, n_{k+1}) + \psi_{k+1}^c(\varepsilon, \varphi, n_{k+1}), \tag{1.22}$$

the coloured lines from the left diagram are slightly perturbed by order of $\varepsilon$ and we can eventually pass over to $\Psi^c(\varepsilon, \varphi, n_{k+1})$ by

$$\Psi^c(\varepsilon, \varphi, n_{k+1}) = \varepsilon^k \cdot \psi_1^c(\varepsilon, \varphi, n_{k+1}) + \cdots + \varepsilon^1 \cdot \psi_k^c(\varepsilon, \varphi, n_{k+1}) + \psi_{k+1}^c(\varepsilon, \varphi, n_{k+1}), \tag{1.23}$$

applying different scalings by order of $\varepsilon^k, \cdots, \varepsilon^1$ to the maps with target subspaces $N_1^c, \cdots, N_k^c$ respectively, as indicated in the middle diagram of figure 1 by the ellipsoid. In this step, the maps with target subspaces $N_{k+1}^c$ and $N_{k+1}$ are not scaled. The origin remains fixed, whereas tangentially touching of the red bended line with the $N_{k+1}$ axis is destroyed by order of $\varepsilon$ during the passage from the left blow up map to the middle diagram with $\varepsilon \neq 0$.

Now, we can turn to the final diffeomorphisms $z = z(\varepsilon) + \varepsilon^k \cdot p_k(\varepsilon) \cdot \{\Psi^c(\varepsilon, \varphi, n_{k+1}) + n_{k+1}\}$ in (1.17) for fixed $\varepsilon \neq 0$. First, the linear near identity transformation $p_k(\varepsilon) = I_B + \varepsilon \cdot \phi_1 + \cdots + \varepsilon^k \cdot \phi_k$ is applied to the middle diagram, again slightly perturbing coloured lines by order of $\varepsilon$ and origin remaining fixed. Secondly, another scaling by $\varepsilon^k$ of all directions is performed and thirdly, the arising ellipsoid type neighbourhood of $0 \in B$ is shifted to the curve position $z(\varepsilon)$ in $B$, as indicated on the right hand side of figure 1.

Reversely, when moving along the curve $z(\varepsilon)$, the ellipsoids are shrinking to zero during the passage from $\varepsilon \neq 0$ to $\varepsilon = 0$ and the blow up system from the left diagram in figure 1 is gradually developing within the ellipsoids in a microscopic version.

Note also that by linearization (1.14), coloured lines in the diagram on the right hand side are representing level sets of $G[\cdot]$, each of which parametrized by $N_{k+1}$.

The linearization (1.14) will be obtained by an application of standard implicit function theorem to a scaled system, ensuring identity of level sets within overlapping images of the $\varepsilon$-parametric family of diffeomorphisms (1.17), as indicated by the grey region in the diagram on the right hand side of figure 1. The union of these images is defining a cone like open neighbourhood of $z(\varepsilon)$ for $\varepsilon \neq 0$ with slowest blow up by $\varepsilon^{2k}$ in the direction of $N_1^c$ and fastest blow up by order of $\varepsilon^k$ in the direction of $N_{k+1}^c \oplus N_{k+1}$.

In case of $k$ =0, no shrinking of the ellipsoids occurs for $\varepsilon \to 0$ and the generalized cone remains a cylinder with linearization and level sets existing in a solid neighbourhood of the origin.

If $N_{k+1} = \{0\}$, bijectivity between $B$ and $\bar{B}$ occurs by (1.10) and each level set degenerates to a single point. In case of $z(\varepsilon) \equiv 0$, $k$-surjectivity occurs with $k = 0$, implying again the constellation of the standard implicit function theorem.

Theorem 1 is valid for general Banach spaces with closed subspaces in (1.10). If a first operator $S_i$ occurs within the sequence of operators $S_1, \cdots, S_i, \cdots, S_{k+1}$, which is of Fredholm type, e.g. $S_1 = G'[0]$ fredholm, then the index of $S_i$ agrees with the dimension of $N_{k+1}$.

In figure 2, two typical generalized cones are qualitatively indicated in $\mathbb{R}^n$, each enclosing its center curve $z(\varepsilon) = \varepsilon^l \cdot z_l + \cdots$, $z_l \neq 0$ of order $l \geq 1$ with $l$ uneven (left) and $l$ even (right).

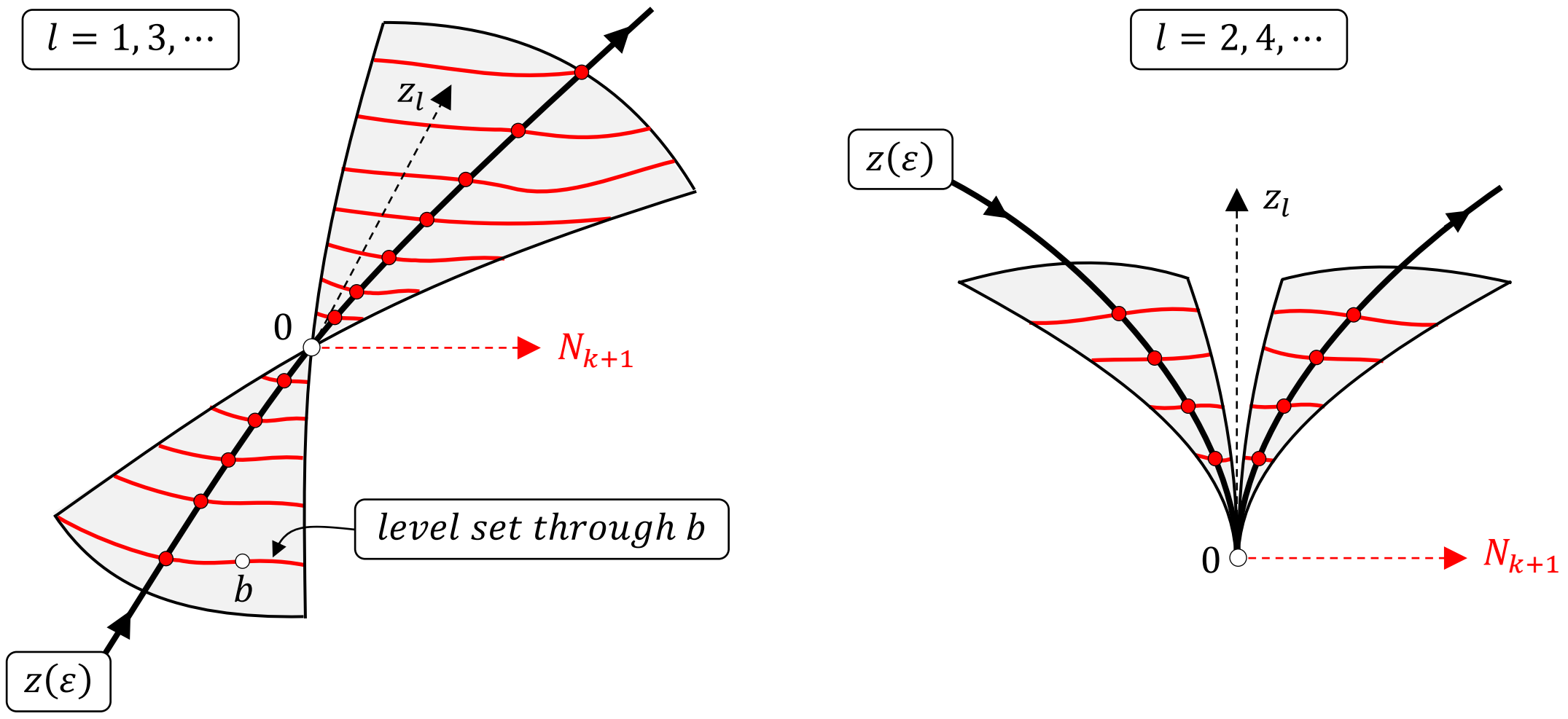

Figure 2: Typical cones with shrinking of subspaces by orders between $\varepsilon^k$ and $\varepsilon^{2k}$.

Geometrically, in both cases red marked level sets are created by moving the nullspace $N_{k+1}$ in a parallel way along the leading direction $z_l \neq 0$ of the curve $z(\varepsilon)$, thereby intersecting the cone and generating the level sets by an $\varepsilon$-perturbation of $N_{k+1}$. In this way, we obtain the level set for every $b$ chosen from the interior of the cone.

On contrary, the behaviour of the level set through $z(0) = 0$ remains, under the previous assumptions, unknown. However, it is precisely this level set that is of interest to us for solving the equation $G[z] = 0$, and thus, for obtaining a generalization of the standard implicit function theorem.

For clarifying the behaviour of the level set through 0, we now aim to derive complete linearization, as stated in (1.6) for the case of $k = 0$. Remember, for $k = 0$, the step from partial linearization (1.4) to complete linearization (1.6) was performed by adding the smallness requirement $\|G[0]\| \ll 1$ to the surjectivity condition $S_1 = G'[0] \in GL[N_1^c, \bar{B}]$.

Now, for general $k \geq 0$, the surjectivity condition $(S_1 \cdots S_{k+1}) \in GL[N_1^c \times \cdots \times N_{k+1}^c, \bar{B}]$ from Theorem 1 is supplemented by the smallness requirement of $z(\varepsilon)$ to be an approximation of order $2k$ for the equation $G[z] = 0$. In detail, we require

$$G[\, z(\varepsilon)\, ] = \varepsilon^{2k} \cdot \bar{b}(\varepsilon) \qquad \textit{with} \qquad \|\, \bar{b}(0)\, \| \ll 1\,. \tag{1.24}$$

Note that for $k = 0$, the first condition is trivially satisfied and the second condition turns into $\|G[0]\| \ll 1$.

A sufficient condition for $\|\bar{b}(0)\| \ll 1$ is given by $\bar{b}(0) = 0$, implying $G[z(\varepsilon)] = \varepsilon^{2k+1} \cdot r(\varepsilon)$ and $z(\varepsilon)$ turns into an approximation of order $2k + 1$. In this sense, condition (1.24) requires $z(\varepsilon)$ to be almost an approximation of order $2k + 1$.

Using (1.24), we can now state a simple Corollary of Theorem 1, concerning complete linearization, as well as a first version of Newton's Lemma in Banach spaces.

**Corollary 1:** Assume the family $L(\varepsilon)$ to be $k$-surjective with all subspaces closed in (1.10).

(i) If the approximation condition (1.24) is satisfied, then the first summand $G[z(\varepsilon)]$ in (1.14) can be absorbed into the linear part by a shift in $\varphi \in N^c$ according to

$$\varphi_0(\varepsilon) := -(S_1 \cdots S_{k+1})^{-1} \cdot \bar{b}(\varepsilon) \in N^c \tag{1.25}$$

and

$$G[\, z(\varepsilon) + \varepsilon^k \cdot p_k(\varepsilon) \cdot \{ \Psi^c(\, \varepsilon,\ \varphi + \varphi_0(\varepsilon),\ n_{k+1}\,) + n_{k+1} \,\}\,]$$

$$= \varepsilon^{2k} \cdot (\, S_1 \cdots S_{k+1}\,) \cdot \varphi\,. \tag{1.26}$$

In the generalized cone

$$z = z(\varepsilon) + \varepsilon^k \cdot p_k(\varepsilon) \cdot \{ \Psi^c(\, \varepsilon, \bar{\varphi}, n_{k+1}\,) + n_{k+1} \,\}\,,\ \ \varepsilon \in U \setminus \{0\}\,,\ \ (\,\bar{\varphi}, n_{k+1}\,) \in V, \tag{1.27}$$

all solutions of $G[z] = 0$ are given by setting $\varphi = 0$ in (1.26) according to

$$z = z(\varepsilon) + \varepsilon^k \cdot p_k(\varepsilon) \cdot \{ \Psi^c(\, \varepsilon,\ \varphi_0(\varepsilon),\ n_{k+1}\,) + n_{k+1} \,\}\,. \tag{1.28}$$

(ii) Assume $k \geq 1$. If $z(\varepsilon)$ is not an approximation of order $2k$ by

$$G[\, z(\varepsilon)\,] = \varepsilon^q \cdot \bar{b}(\varepsilon) \quad \textit{with} \quad 0 \leq q \leq 2k-1 \quad \textit{and} \quad \bar{b}(0) \neq 0\,, \tag{1.29}$$

then no zero of $G[z]$ occurs in the cone (1.27).

From a geometrical viewpoint, low order approximation of $z(\varepsilon)$ according to (1.29) implies the subspace $N_{k+1}$ and the generalized cone (1.27) having no common elements, as depicted in figure 2 above. Then by (ii), a possibly existing level set through the origin is not contained in the cone.

However, when successively increasing the order of approximation from $\varepsilon^q, 0 \leq q \leq 2k-1$ up to $\varepsilon^{2k}$ as required in (i), then the subspace $N_{k+1}$ is gradually turning into the direction of the cone until condition (1.24) forces $N_{k+1}$ to point into the interior of the cone, as qualitatively shown in figure 3.

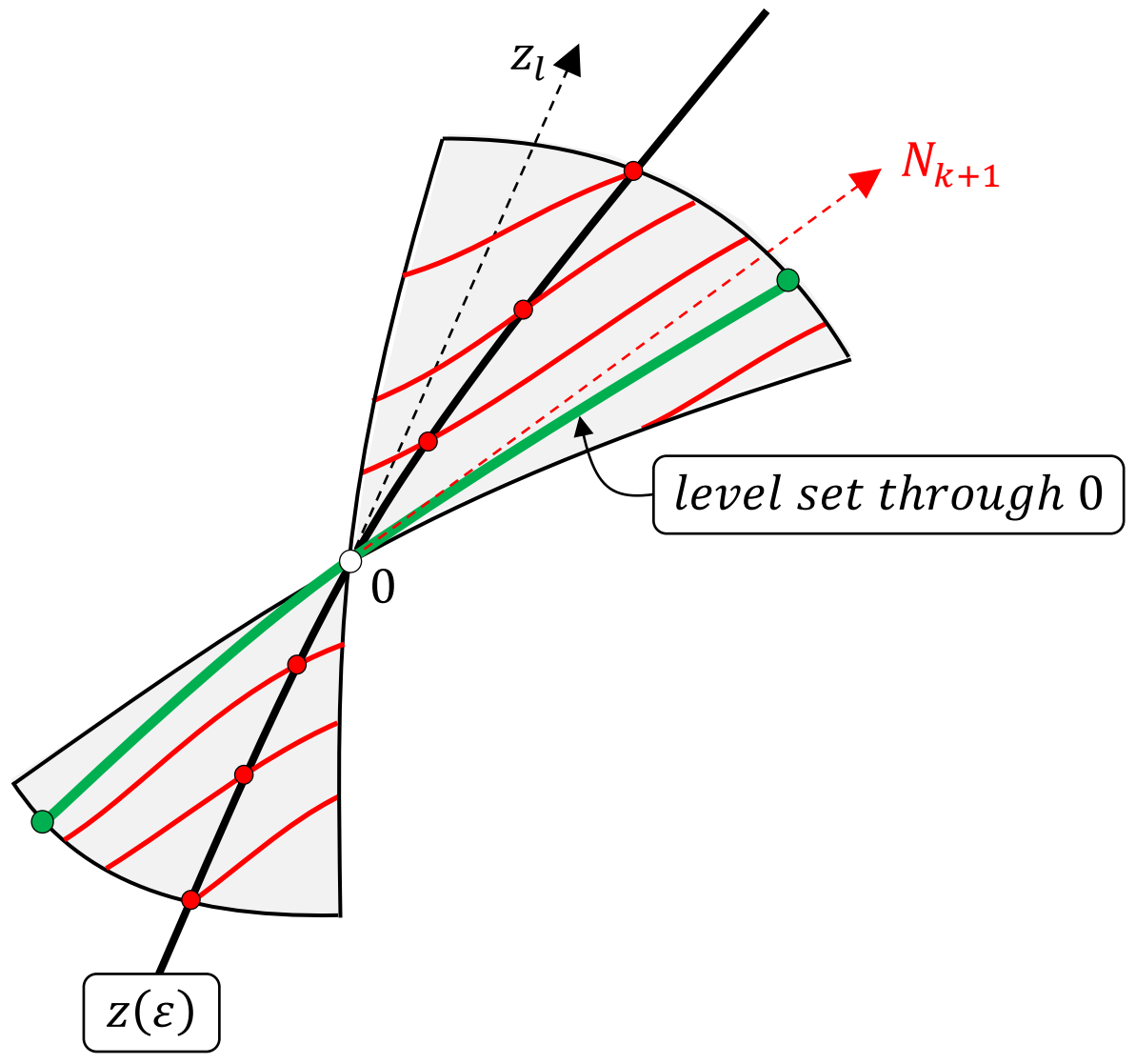

Figure 3: Green marked level set through the origin in case of $N_{k+1}$ pointing into the cone.

Then, by moving $N_{k+1}$ in a parallel way in the direction of $z_l$, the level sets are again obtained by parametrization with respect to $N_{k+1}$, but now including the green marked level set through the origin.

The solutions (1.28) of $G[z]=0$ in Corollary 1 are derived after complete linearization is achieved by construction of the transformation $p_k(\varepsilon)\cdot\{\Psi^c(\varepsilon,\varphi,n_{k+1})+n_{k+1}\}$ in Theorem 1. However, by a direct analysis of the equation $G[z]=0$, i.e. without prior transformation, the following improved version of Newton's Lemma can be obtained.

**Corollary 2:** Assume $z(\varepsilon)$ to be an approximation of order $2k$ by $G[z(\varepsilon)]=\varepsilon^{2k}\cdot\bar{b}(\varepsilon)$ and the family $L(\varepsilon)$ to be $k$-surjective with $k\geq 1$ and all subspaces closed in (1.10).

Let $\mathcal{P}_i$ denote the projection to $R_i, i=1,\cdots,k+1$ with respect to the direct sum of $\bar{B}$ in (1.10) and consider the first bilinear term $\mathcal{B}_0 := \frac{1}{2}G''[0]$ within expansion (1.1).

(i) Assume smallness of the remainder $\bar{b}(\varepsilon)$ with respect to the subspace $R_{k+1}$ according to

$$\|\,\mathcal{P}_{k+1}\cdot\bar{b}(0)\,\| \ll 1\,. \tag{1.30}$$

Then, for $\varepsilon\in U$ and $\|\,n_{k+1}\,\|\ll 1$, a smooth mapping of the form

$$\hat{n}^c(\varepsilon,n_{k+1}) = \varepsilon^k\cdot\underbrace{\bar{n}_1^c(\varepsilon,n_{k+1})}_{\to N_1^c}+\cdots+\varepsilon^1\cdot\underbrace{\bar{n}_k^c(\varepsilon,n_{k+1})}_{\to N_k^c}+\underbrace{\bar{n}_{k+1}^c(\varepsilon,n_{k+1})}_{\to N_{k+1}^c},$$

exists, such that

$$z = z(\varepsilon)+\varepsilon^k\cdot p_k(\varepsilon)\cdot\{\,\hat{n}^c(\,\varepsilon,n_{k+1}\,)+n_{k+1}\,\} \tag{1.31}$$

defines solutions of $G[z]=0$.

For $\varepsilon=0$, the last summand $\bar{n}_{k+1}^c(0,n_{k+1})\in N_{k+1}^c$ of the mapping $\hat{n}^c(\varepsilon,n_{k+1})$ represents an existing solution of the quadratic equation

$$S_{k+1}\cdot n_{k+1}^c+\mathcal{P}_{k+1}\cdot\mathcal{B}_0\cdot(\,n_{k+1}^c+n_{k+1}\,)^2+\mathcal{P}_{k+1}\cdot\bar{b}(0)\;=0\,, \tag{1.32}$$

satisfying $\|\,\bar{n}_{k+1}^c(\varepsilon,n_{k+1})\,\|\ll 1$, whereas first $k$ summands read

$$\begin{aligned}\bar{n}_1^c(0,n_{k+1}) &= -S_1^{-1}\cdot\{\,\mathcal{P}_1\cdot\mathcal{B}_0\cdot[\,\bar{n}_{k+1}^c(0,n_{k+1})+n_{k+1}\,]^2+\mathcal{P}_1\cdot\bar{b}(0)\,\}\in N_1^c\\ &\;\;\vdots\\ \bar{n}_k^c(0,n_{k+1}) &= -S_k^{-1}\cdot\{\,\mathcal{P}_k\cdot\mathcal{B}_0\cdot[\,\bar{n}_{k+1}^c(0,n_{k+1})+n_{k+1}\,]^2+\mathcal{P}_k\cdot\bar{b}(0)\,\}\in N_k^c\,.\end{aligned} \tag{1.33}$$

(ii) Assume the curvature of $G$ at $z=0$ in the direction of $N_{k+1}^c$ equals zero according to

$$\mathcal{B}_0\cdot(\,n_{k+1}^c\,)^2 = 0\,. \tag{1.34}$$

Then, for $\varepsilon\in U$ and $\|\,n_{k+1}\,\|\ll 1$, a smooth mapping of the form

$$\tilde{n}^c(\varepsilon,n_{k+1}) = \varepsilon^k\cdot\underbrace{\bar{n}_1^c(\varepsilon,n_{k+1})}_{\to N_1^c}+\cdots+\varepsilon^1\cdot\underbrace{\bar{n}_k^c(\varepsilon,n_{k+1})}_{\to N_k^c}+\underbrace{\bar{n}_{k+1}^c(\varepsilon,n_{k+1})}_{\to N_{k+1}^c},$$

exists, such that

$$z = z(\varepsilon)+\varepsilon^k\cdot p_k(\varepsilon)\cdot\{\,\tilde{n}^c(\,\varepsilon,n_{k+1}\,)+n_{k+1}\,\} \tag{1.35}$$

defines solutions of $G[z]=0$.

For $\varepsilon = 0$, the last summand $\bar{n}^c_{k+1}(0, n_{k+1}) \in N^c_{k+1}$ of the mapping $\tilde{n}^c(\varepsilon, n_{k+1})$ is explicitly given by

$$\bar{n}^c_{k+1}(0, n_{k+1}) = -(\,S_{k+1} + 2\cdot\mathcal{P}_{k+1}\cdot\mathcal{B}_0\cdot n_{k+1}\,)^{-1}\cdot\mathcal{P}_{k+1}\cdot\left[\,\mathcal{B}_0\cdot(n_{k+1})^2 + \bar{b}(0)\,\right], \quad (1.36)$$

whereas first $k$ summands are again determined by (1.33).

Compared to the smallness requirement (1.24) of Corollary 1, within the improved version of Corollary 2 (i), smallness of the remainder is required only with respect to the subspace $R_{k+1}$, i.e. the remainder

$$\bar{b}(0) \;=\; \underbrace{\mathcal{P}_1\cdot\bar{b}(0) + \cdots + \mathcal{P}_k\cdot\bar{b}(0)}_{arbitrary} \;+\; \underbrace{\mathcal{P}_{k+1}\cdot\bar{b}(0)}_{small} \quad (1.37)$$

is allowed to be an arbitrary element in $\bar{B}$, with the restriction that it is almost positioned within $R_1 \oplus \cdots \oplus R_k$. As a consequence, first $k$ components $\bar{n}^c_1(\varepsilon, n_{k+1}), \cdots, \bar{n}^c_k(\varepsilon, n_{k+1})$ of the solutions are now allowed to take arbitrary values, caused by arbitrary values of $-S_1^{-1}\cdot\mathcal{P}_1\cdot\bar{b}(0), \cdots,$ $-S_k^{-1}\cdot\mathcal{P}_k\cdot\bar{b}(0)$ in (1.33).

In Theorem 1, Corollary 1 and Corollary 2 (i), the results are obtained by considering certain properties of $G[z(\varepsilon)]$ and the linear term $L(\varepsilon)\cdot b$ within the expansion (1.1). These assumptions are chosen sufficiently strong for complete control of higher order terms $\mathcal{B}(\varepsilon)\cdot b^2 + r(\varepsilon, b)\cdot b^3$. Now, if the impact of higher order terms on the action of the operator $G[z]$ is decreased by appropriate simplifications, then it should also be possible to weaken the assumptions concerning the leading terms $G[z(\varepsilon)]$ and $L(\varepsilon)$.

In this spirit, Corollary 2 (ii) assumes the leading nonlinear term $\mathcal{B}_0 = \frac{1}{2}G''[0]$ in (1.1) to be simplified in such a way that the curvature of $G$ at $z = 0$, taken with respect to the direction of $N^c_{k+1} \subset B$, equals zero, i.e. the impact of $\mathcal{B}_0$ on the action of $G[z]$ is weakened according to

$$\underbrace{\mathcal{B}_0\cdot b^2}_{curvature} \;=\; \mathcal{B}_0\cdot(n^c_1 + \cdots + n^c_{k+1} + n_{k+1})^2 \;=\; \underbrace{\mathcal{B}_0\cdot(n^c_{k+1})^2}_{=0} \;+\; \underbrace{\mathcal{B}_0\cdot[\,(n^c_1)^2 + \cdots\,]}_{arbitrary}\,. \quad (1.38)$$

Then, by Corollary 2 (ii), we see that by assuming (1.38), it is no longer necessary to require the projection of $G[z(\varepsilon)] = \varepsilon^{2k}\cdot\bar{b}(\varepsilon)$ to $R_{k+1}$ to be small according to (1.30), thus weakening the approximation requirement for $z(\varepsilon)$ to an arbitrary approximation of order $\varepsilon^{2k}$.

In section 6 of this paper, the impact of higher order terms in the expansion (1.1) will systematically be reduced, yielding the possibility of gradually weakening the approximation requirements of $G[z(\varepsilon)]$ necessary to obtain solutions of $G[z] = 0$. Along these lines, needed approximation can be relaxed from approximations of order $\varepsilon^{2k}$ down to approximations of order $\varepsilon^{k+1}$. The size of the cone increases correspondingly.

Some of the results in Theorem 1 and Corollary 1, 2 are given in [15] and [17], however, without explicit use of Jordan chains. The consistent derivation stressing Jordan chains in the paper at hand, might be new, motivated by [2], [5], [11], [12] and [7].

Results concerning approximations $z(\varepsilon)$ of order $2k+1$ by $G[z(\varepsilon)] = \varepsilon^{2k+1}\cdot\bar{b}(\varepsilon)$, combined with nondegeneracy of order $k$, occur in several places of algebra and analysis, applied to mappings of quite different properties with respect to smoothness and type of participating spaces.

Compare Newton's Lemma [8], [3], Tougeron's implicit function theorem [4], [9], [14] or Hensel's Lemma [6].

In this paper, the constellation is applied to $G[z] = 0,\ G \in C^r(B, \bar{B})$ with $B, \bar{B}$ real or complex Banach spaces. In a first step, we are weakening needed approximation from $G[z(\varepsilon)] = \varepsilon^{2k+1} \cdot \bar{b}(\varepsilon)$ down to $G[z(\varepsilon)] = \varepsilon^{2k} \cdot \bar{b}(\varepsilon)$ and $\|\mathcal{P}_{k+1} \cdot \bar{b}(0)\| \ll 1$, as stated in Corollary 2 (i). This results from a detailed analysis of the blown up remainder equation under consideration of the first bilinear term $\mathcal{B}_0 \cdot b^2 = \frac{1}{2} G''[0] \cdot b^2$ in (1.1). Note also that by (1.28), (1.31) and (1.35), given approximation $z(\varepsilon)$ and resulting solutions $z(\varepsilon) + \varepsilon^k \cdot p_k(\varepsilon) \cdot \{\Psi^c(\,\cdot\,) + n_{k+1}\}$ merely have to agree up to the order of $k-1$, in contrast to standard Newton Lemma with identity up to the order of $k$. The cases $\mathbb{K} = \mathbb{R}$ and $\mathbb{K} = \mathbb{C}$ are treated in a parallel way.

In section 2, results concerning Jordan chains with rank between 0 and $k$ are collected from [11] and [18], used to define a generalized cone around the curve $z(\varepsilon)$, which is composed of approximations of order $2k$, all of which possessing surjectivity of order $k$.

Based on this homogenous cone, Theorem 1 is proved in section 3 by a blow up transformation and subsequent application of the standard implicit function theorem to a regular remainder equation. Corollary 1, i.e. our base version of Newton's Lemma, follows immediately in section 4.

Sections 5 and 6 are dealing more intensively with the Lemma of Newton. In section 5, Corollary 2 is proved, whereas in section 6, further refinements of Newton's Lemma are derived in a systematic way, as indicated above.

At the end of section 6, we use the example $G[x, y] = -xy^3 + x^5 + y^5$ of an extended $E_6$ singularity [1] for calculation of $k$-surjectivity and local zero sets. Here, two solution curves emanate from the origin possessing order of surjectivity given by $k = 3$ and $k = 11$. Using this example, also the relation to a constant topological degree in each of the half cones is established and the Milnor number $\mu$ of the singularity is calculated from $k$-values of the two solution curves. For examples considering nonlinear systems of differential algebraic equations, we refer to [19].

In section 7, the redundancy of charts, present within the $\varepsilon$-parametric family of diffeomorphisms (1.17), is eliminated by an appropriate parametrization condition, guided by the leading order term of the center curve $z(\varepsilon) = \varepsilon^l \cdot z_l + \cdots$ satisfying $z_l \neq 0$. Then, each level set in the cone will be uniquely parametrized with respect to $\varepsilon \in U$ and a subspace $\Pi_{k+1} \subsetneq N_{k+1}$ of codimension 1.

We close section 1 with two remarks concerning perturbation and homogeneity of the cones.

**Remarks 1)** How do the results of Theorem 1 and Corollaries 1, 2 change under perturbation of the operator $G[z]$ or the curve $z(\varepsilon)$? The answer is easily obtained from the proofs, by looking for maximal derivatives of $G[z]$ and $z(\varepsilon)$, needed to establish a certain property. Then, perturbing above these maximal derivatives, will not change the property.

Without going into details, we state the results concerning perturbation of $G[z]$ and $z(\varepsilon)$ in the following Lemma; compare also [16], [17].

**Lemma 1:** Assume the family $L(\varepsilon)$ to be $k$-surjective with $k \geq 1$ and all subspaces closed within the decompositions (1.10).

(i) Perturbing the mapping $G[z]$ by order of $i+1,\ i \in \{0, \cdots, k\}$ according to

$$\bar{G}[z] \coloneqq G[z] + \alpha \cdot \left( M_{i+1} \cdot z^{i+1} + \cdots + M_{k+1} \cdot z^{k+1} \right) + O\left( z^{k+2} \right), \tag{1.39}$$

$\alpha \in \mathbb{K}$ chosen sufficiently small, implies the order $\bar{k}$ of surjectivity of the perturbed linearized family $\bar{L}(\varepsilon) := \bar{G}'[z(\varepsilon)]$ to satisfy $i \le \bar{k} \le k$. In case of $\alpha = 0$, the perturbation reads $\bar{G}[z] = G[z] + O\big(z^{k+2}\big)$ and the family $\bar{L}(\varepsilon)$ preserves $k$-surjectivity of $L(\varepsilon)$.

The mappings $M_\tau \in L^\tau[B, \bar{B}]$ denote $\tau$-linear mappings and $O\big(z^{k+2}\big)$ denotes a remainder function $r(z)$ satisfying $\|r(z)\| < c \cdot \|z\|^{k+2}$ with $c > 0$.

(ii) Perturbing the curve $z(\varepsilon)$ by order of $i$, $i \in \{1, \cdots, k\}$ according to

$$\bar{z}(\varepsilon) := z(\varepsilon) + \alpha \cdot \big(\, \varepsilon^i \cdot b_i + \cdots + \varepsilon^k \cdot b_k \,\big) + O\big(\, \varepsilon^{k+1} \,\big)\,, \tag{1.40}$$

$\alpha \in \mathbb{K}$ chosen sufficiently small, implies the order $\bar{k}$ of surjectivity of the perturbed linearized family $\bar{L}(\varepsilon) := G'[\bar{z}(\varepsilon)]$ to satisfy $i \le \bar{k} \le k$. In case of $\alpha = 0$, the family $\bar{L}(\varepsilon)$ retains $k$-surjectivity of $L(\varepsilon)$. Again, $O(\varepsilon^{k+1})$ denotes a remainder function $r(\varepsilon)$ satisfying $\|r(\varepsilon)\| < c \cdot |\varepsilon|^{k+1}$ with $c > 0$.

(iii) Assume the approximation condition $G[z(\varepsilon)] = \varepsilon^{2k} \cdot \bar{b}(\varepsilon)$ with $\|\bar{b}(0)\| \ll 1$ from (1.24) to be satisfied. Then, perturbing $G[z]$ and $z(\varepsilon)$ according to

$$\bar{G}[z] := G[z] + \alpha_1 \cdot M_{2k} \cdot z^{2k} \;+\; O\big(\, z^{2k+1} \,\big) \tag{1.41}$$

$$\bar{z}(\varepsilon) := z(\varepsilon) + \alpha_2 \cdot \varepsilon^{2k} \cdot b_{2k} \;+\; O\big(\, \varepsilon^{2k+1} \,\big)\,, \tag{1.42}$$

$\alpha_1, \alpha_2 \in \mathbb{K}$ chosen sufficiently small, implies (1.24) to remain true with $\bar{G}[\bar{z}(\varepsilon)] = \varepsilon^{2k} \cdot \tilde{b}(\varepsilon)$ and $\|\tilde{b}(0)\| \ll 1$.

In case of (i), perturbation of $G[z]$ with $\alpha = 0$, does not effect $k$-surjectivity of $z(\varepsilon)$ for arbitrary magnitude of the perturbation $r(z) = O\big(z^{k+2}\big)$. Also, perturbation of $G[z]$ by order below $z^{k+2}$ cannot destroy the property of surjectivity for $\alpha$ chosen sufficiently small, however, typically lowering the order at which surjectivity of $\bar{L}(\varepsilon)$ occurs, down to a value of $\bar{k}$ between $i$ and $k$. The magnitude of the perturbation has to be small, which is controlled by $\alpha$ chosen sufficiently small.

Basically, the same behaviour takes place for mapping $G[z]$ not perturbed, but instead curve $z(\varepsilon)$ perturbed, as stated in (1.40). Perturbation of $z(\varepsilon)$ by $O\big(\varepsilon^{k+1}\big)$ does not effect $k$-surjectivity of the new family $\bar{L}(\varepsilon)$. Perturbation of $z(\varepsilon)$ by order below $\varepsilon^{k+1}$ cannot destroy the property of surjectivity for $\alpha$ chosen sufficiently small, but may decrease the order of surjectivity of $\bar{L}(\varepsilon)$, again to a value of $\bar{k}$ between $i$ and $k$.

Concerning the approximation condition in (iii), small perturbation of order $2k$ and arbitrary perturbation of order $2k+1$ will not destroy (1.24), and hence the solutions of $G[z] = 0$ from Corollary 1 (i) carry over to solutions of $\bar{G}[z] = 0$, only varying higher order derivatives within the solution manifold. Qualitatively, the same holds true, when replacing $\|\bar{b}(0)\| \ll 1$ from (1.24) by the weaker assumption $\|\mathcal{P}_{k+1} \cdot \bar{b}(0)\| \ll 1$ from (1.30) in Corollary 2 (i).

In summary, perturbing $G[z]$ or $z(\varepsilon)$ according to (1.39), (1.40) does not destroy the property of surjectivity of the family and, in particular, the perturbation is not capable to increase the value of $k$. In other words, small perturbation controlled by $\alpha$ means a step towards a more regular constellation, and hence, the order $k$ of surjectivity is showing semicontinous behaviour from above under perturbation (comparable to the algebraic multiplicity of a classical eigenvalue).

**2)** Assume the family $L(\varepsilon) = G'[z(\varepsilon)]$ to be $k$-surjective with $k \geq 1$. Then, the generalized cone

$$z_c = z(\varepsilon) + \varepsilon^k \cdot p_k(\varepsilon) \cdot \{ \Psi^c(\varepsilon, \varphi, n_{k+1}) + n_{k+1} \} \tag{1.43}$$

from Theorem 1 is well defined, enclosing its center line $z(\varepsilon)$. Now, take an arbitrary smooth curve within the cone defined for $\varepsilon \in U$ by $\varphi(\varepsilon)$ and $n_{k+1}(\varepsilon)$ according to

$$z_c(\varepsilon) \coloneqq z(\varepsilon) + \varepsilon^k \cdot p_k(\varepsilon) \cdot \{ \Psi^c [\varepsilon, \varphi(\varepsilon), n_{k+1}(\varepsilon)] + n_{k+1}(\varepsilon) \}. \tag{1.44}$$

The curve $z_c(\varepsilon)$ represents a perturbation of the center line $z(\varepsilon)$ of order $\varepsilon^k$ and it is straightforward to see that the perturbed linearized family $\bar{L}(\varepsilon) = G'[z_c(\varepsilon)]$ along $z_c(\varepsilon)$ adopts order $k$ of surjectivity from the center line $z(\varepsilon)$.

In the upper diagram of figure 4, the constellation is qualitatively indicated in $\mathbb{R}^n$.

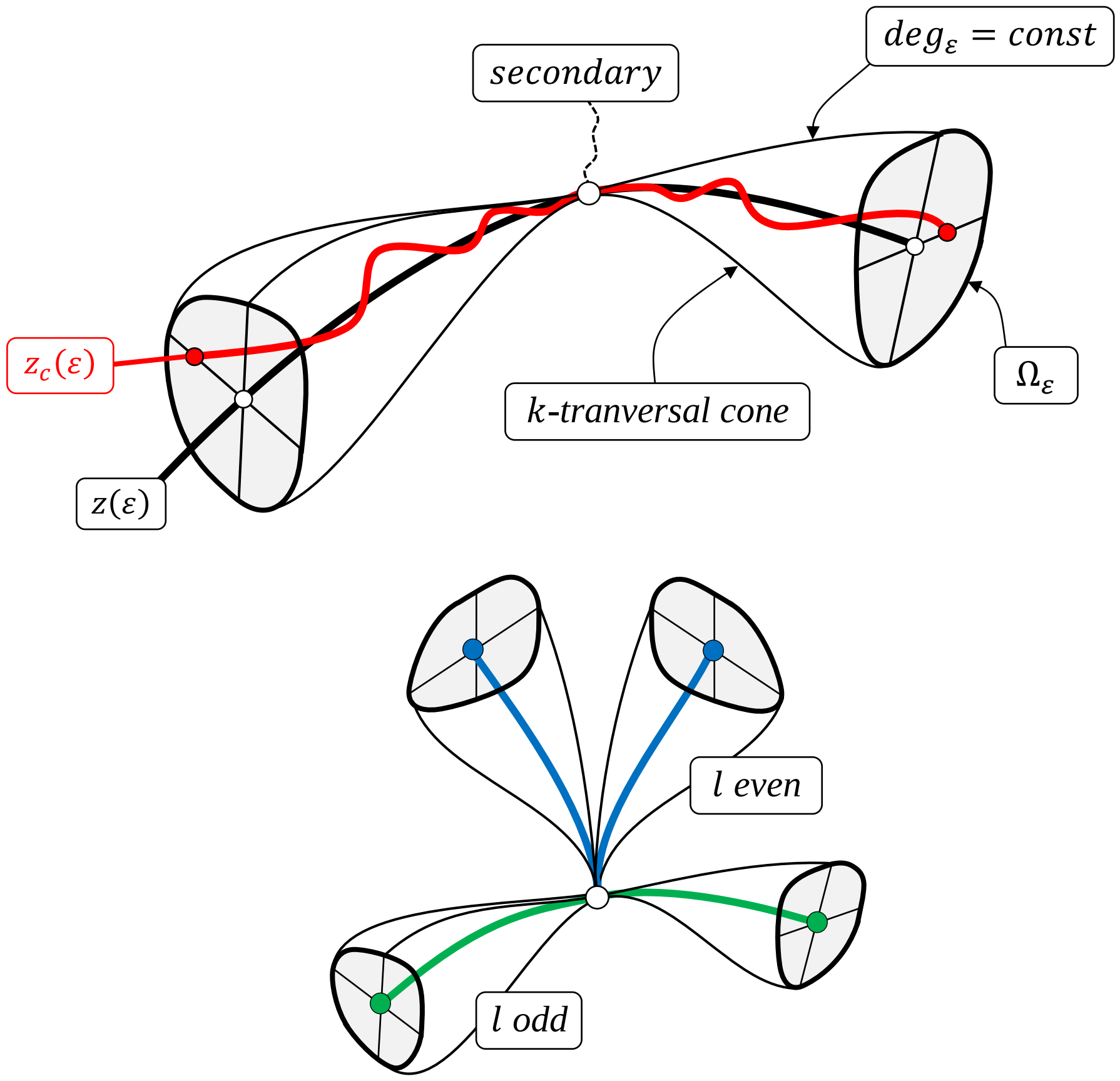

Figure 4: Homogeneity within a $k$-transversal cone and two solution curves near a singularity.

Hence, we see that $k$-surjectivity is not only a property of the center line $z(\varepsilon)$, but a property of the cone by itself, and along these lines, the cone may be denoted as a $k$-transversal cone [17]. In addition, if the center line $z(\varepsilon)$ is an approximation of order $2k$ by $G[z(\varepsilon)] = \varepsilon^{2k} \cdot \bar{b}(\varepsilon)$, then by Theorem 1 (i), all curves $z_c(\varepsilon)$ in the cone given by (1.44) are approximations of order $2k$ as well.

Finally, if $\mathbb{K} = \mathbb{R}$ and $B, \bar{B}$ are of finite dimensions, then the topolocical degree $deg_\varepsilon$ is well defined along appropriate cross sections $\Omega_\varepsilon$ (grey shaded), satisfying $deg_\varepsilon = const$ in each of the half cones, if $dim\ N_{k+1} = 1$ and $G[z(\varepsilon)] = \varepsilon^{2k+1} \cdot \bar{b}(\varepsilon)$ are presupposed [17].

Using this information, the two half cones are perfectly suited for investigation of possible sign change of the topological degree and subsequent bifurcation of secondary solutions outside the cone, as shown in [17]. Far reaching local and global bifurcation results, using the topological degree, can be found in [11], [12] in the context of Fredholm operators and with respect to mappings of the form $G[\lambda, y] = L(\lambda) \cdot y + O(y^2)$. In the upper diagram of figure 4, a dashed line indicates some secondary solutions outside the cone.

In the lower diagram of figure 4, a typical constellation near a singularity is indicated with blue and green marked solution curves of $G[z] = 0$, each enclosed by half cones possessing constant topological degree. The order $k$ of surjectivity will be different for the two curves, where the leading coefficient $z_l$ of the cusp type blue curve is characterized by an even number $l = 2, 4, \dots$ , whereas the green curve corresponds to an odd number $l = 1, 3, \dots$ .

Summing up, a generalized cone is characterized by homogenous or trivial behaviour in the sense that level sets can be straightened out to obtain linearization, and secondly, for every curve given by (1.44), $k$-surjectivity and approximation of order $2k$ can be ensured. Moreover, each of the half cones may possess a constant topological degree, possibly usable to ensure secondary bifurcation; compare also the example at the end of section 6.

## *2. Jordan chains*

The linearization Theorem 1 relies on precise knowledge about the behaviour of the family of linear operators $L(\varepsilon) \in L[B, \bar{B}]$, yielding control of the higher order terms $\mathcal{B}(\varepsilon) \cdot b^2 + r(\varepsilon, b) \cdot b^3$ within the expansion (1.1). If this control is sufficiently strong, ensured by $k$-surjectivity of $L(\varepsilon)$, then higher order terms can be eliminated by pure right transformation given by (1.17).

Geometrically, $\bar{b} = L(\varepsilon) \cdot b$ gives the directional derivative of $G[z]$ at $z = z(\varepsilon)$ in the direction of $b$. At $\varepsilon = 0$, the directional derivative mapping $L(0) = L_0 = G'[0]$ is allowed to be highly degenerate, e.g. $L_0 \equiv 0$, however turning into a surjective operator for $\varepsilon \neq 0$ and it remains to understand in detail, how surjectivity arises during the passage from $\varepsilon = 0$ to $\varepsilon \neq 0$.

Now, in finite dimensions (or with Fredholm operators and a Lyapunov-Schmidt reduction), the main tool for analyzing this passage is given by the behaviour of some determinant of the linearization $L(\varepsilon) = G'[z(\varepsilon)]$, constructed with respect to a transversal subspace $N^c$ of $z(\varepsilon)$, i.e.

$$det\{ G_{N^c}[ z(\varepsilon) ] \} = \varepsilon^{\chi} \cdot r(\varepsilon) \quad with \quad \chi \geq 0 \quad and \quad r(0) \neq 0 \, . \tag{2.1}$$

Here $G_{N^c}[\cdot]$ denotes linearization in the direction of $N^c$. Note that the determinant measures the overall speed of expansion, when mapping the unit ball from $N^c \subset B$ to $\bar{B}$, i.e. in (2.1) the expansion is characterized by $\varepsilon^{\chi}$. Hence, a suitable Ansatz in finite dimensions for proving Newton's Lemma [8], Hensel's Lemma [6] or Tougeron's implicit function theorem [9] is given by

$$z \ = \ z(\varepsilon) \ + \ det\{ G_{N^c}[ z(\varepsilon) ] \} \cdot \bar{z} \ = \ z(\varepsilon) \ + \ \varepsilon^{\chi} \cdot r(\varepsilon) \cdot \bar{z} \, , \tag{2.2}$$

as well as factoring out of the adjoint matrix of the Jacobian $G_{N^c}[z(\varepsilon)]$. Then, a desingularized remainder equation occurs, appropriate for application of standard implicit function theorem with respect to the new variable $\bar{z}$. Note also that condition $r(0) \neq 0$ in (2.1) represents a nondegeneracy condition, ascertaining an isolated singularity to occur at $\varepsilon = 0$.

In infinite dimensions, the determinant and the adjoint matrix are no longer present in a straightforward way, implying the nondegeneracy condition and the proof by itself to be refor-mulated appropriately. Then, the main tool is given by curves $b(\varepsilon) = b_0 + \varepsilon \cdot b_1 + \cdots$ in $B$ that are mapped to $\bar{B}$ by the linearization $G'[z(\varepsilon)] = L(\varepsilon) = L_0 + \varepsilon \cdot L_1 + \cdots$ according to

$$L(\varepsilon) \cdot b(\varepsilon) = \underbrace{L_0 b_0}_{0-th\ order} + \underbrace{\varepsilon \cdot (L_0\ L_1) \begin{pmatrix} b_1 \\ b_0 \end{pmatrix}}_{1-th\ order} + \underbrace{\varepsilon^2 \cdot (L_0\ L_1\ L_2) \begin{pmatrix} b_2 \\ b_1 \\ b_0 \end{pmatrix}}_{2-th\ order} + \cdots . \tag{2.3}$$

If $L_0 = G'[0]$ is surjective, then the 0-th order coefficient $L_0 b_0$ takes every element in $\bar{B}$, implying 0-surjectivity of the family $L(\varepsilon)$.

If $L_0$ is not surjective by $R[L_0] \subsetneq \bar{B}$, then we may add leading coefficients of 1-th order given by $\bar{R}_1 := \{L_0 b_1 + L_1 b_0 \mid b_1 \in B,\ L_0 b_0 = 0\}$ and, if $\bar{R}_1 = \bar{B}$, we obtain 1-surjectivity of $L(\varepsilon)$.

Now, in general a leading coefficient of order $k \geq 1$ results from a curve $b(\varepsilon)$, which is an approximation of order $k$ with respect to the equation $L(\varepsilon) \cdot b = 0$, i.e. a $k$-th order leading coefficient arises from a constellation of the following form

$$L(\varepsilon) \cdot b(\varepsilon) = \underbrace{\sum_{l=0}^{k-1} \varepsilon^l \cdot \sum_{i+j=l} L_i b_j}_{=0} + \underbrace{\varepsilon^k \cdot (L_0\ L_1 \cdots\ L_k) \cdot \begin{pmatrix} b_k \\ b_{k-1} \\ \vdots \\ b_0 \end{pmatrix}}_{\substack{k-th\ order \\ leading\ coefficients}} + \sum_{l=k+1}^{\infty} \varepsilon^l \cdot \sum_{i+j=l} L_i b_j$$

$$=: \varepsilon^k \cdot \bar{b}(\varepsilon) \tag{2.4}$$

and we see that first $k$ summands have to vanish for a $k$-th order leading coefficient to appear.

Red marked coefficients $(b_0, \cdots, b_{k-1})$, $b_0 \neq 0$ of an approximation $b(\varepsilon) = b_0 + \cdots + \varepsilon^{k-1} \cdot b_{k-1} + \cdots$ of order $k \geq 1$ are called a Jordan chain of length $k$ (or chain of generalized eigenvectors of the eigenvalue $\varepsilon = 0$). The base element $b_0 \in N[L_0]$, $b_0 \neq 0$ is called the root element of the Jordan chain.

Finally, the maximal order of approximation that can be constructed out of $b \in B$ is called the rank of $b$ with abbreviation $rk(b)$. Note that $rk(b) = 0$, if $b \notin N[L_0]$ and $rk(0) = \infty$.

Now, from [18] we know that $k$-surjectivity of $L(\varepsilon)$ is equivalent to the existence of direct sums (1.10) and a polynomial family of linear mappings

$$p_k(\varepsilon) = \underbrace{\phi_0}_{= I_B} + \varepsilon \cdot \phi_1 + \cdots + \varepsilon^k \cdot \phi_k\ , \quad \phi_i \in L[\,B, N_1^c \oplus \cdots \oplus N_{k+1-i}^c\,]\ , \quad i = 1, \ldots, k \tag{2.5}$$

such that

$$L(\varepsilon) \cdot p_k(\varepsilon) \cdot \overbrace{n_{i+1}^c}^{rk=i} = \varepsilon^i \cdot S_{i+1} \cdot n_{i+1}^c + O(\,\varepsilon^{i+1}\,)\ , \quad i = 0, \ldots, k \tag{2.6}$$

$$L(\varepsilon) \cdot p_k(\varepsilon) \cdot n_{k+1} = O(\,\varepsilon^{k+1}\,) \tag{2.7}$$

with $S_{i+1} \in GL[N_{i+1}^c, R_{i+1}]$, $i = 0, \ldots, k$ and matrix representation of the operator polynomial $P_k(\varepsilon) \coloneqq S_1 + \cdots + \varepsilon^k \cdot S_{k+1}$ given by lower triangularity according to

$$
\begin{array}{c} \\ P_k(\varepsilon) = \end{array}
\begin{array}{cccc} \boxed{N_1^c} & \cdots & \boxed{N_{k+1}^c} & \boxed{N_{k+1}} \\ \end{array}
$$

$$
P_k(\varepsilon) = \begin{pmatrix} \boxed{S_1} & \square & \square & \square \\ \boxtimes & \boxed{\ddots} & \square & \square \\ \boxtimes & \boxtimes & \boxed{\varepsilon^k S_{k+1}} & \square \end{pmatrix} \begin{array}{c} \boxed{R_1} \\ \vdots \\ \boxed{R_{k+1}} \end{array} . \tag{2.8}
$$

Squares without entry denote the zero operator, red marked entries indicate bijection for $\varepsilon \neq 0$, while black crosses denote possible entry into subspaces of $\bar{B}$, as indicated on the right.

By (2.6), the mapping $b(\varepsilon) = p_k(\varepsilon) \cdot n_{i+1}^c$ represents an approximation of order $i$ to $L(\varepsilon) \cdot b = 0$ with rank satisfying $rk(n_{i+1}^c) = i$ for $n_{i+1}^c \neq 0$ and $i = 0, \ldots, k$. In addition by (2.7), the mapping $b(\varepsilon) = p_k(\varepsilon) \cdot n_{k+1}$ represents an approximation of order $k+1$ (with rank unspecified).

From a geometrical point of view, we see from (2.5), (2.6) that, starting at $\varepsilon = 0$ with $L_0 \cdot n_{i+1}^c = 0$ for $i = 1, \ldots, k$, the root element $n_{i+1}^c$ can smoothly be varied by $p_k(\varepsilon) \cdot n_{i+1}^c$ in such a way that by $L(\varepsilon) \cdot p_k(\varepsilon) \cdot n_{i+1}^c$ a subspace is created in $\bar{B}$, which blows up with $\varepsilon$-speed of order $i$.

In this sense, the Jordan chains describe in a perfect way the passage from nonsurjectivity at $\varepsilon = 0$ to surjectivity for $\varepsilon \neq 0$. The situation is indicated in figure 5 within $\mathbb{R}^n$.

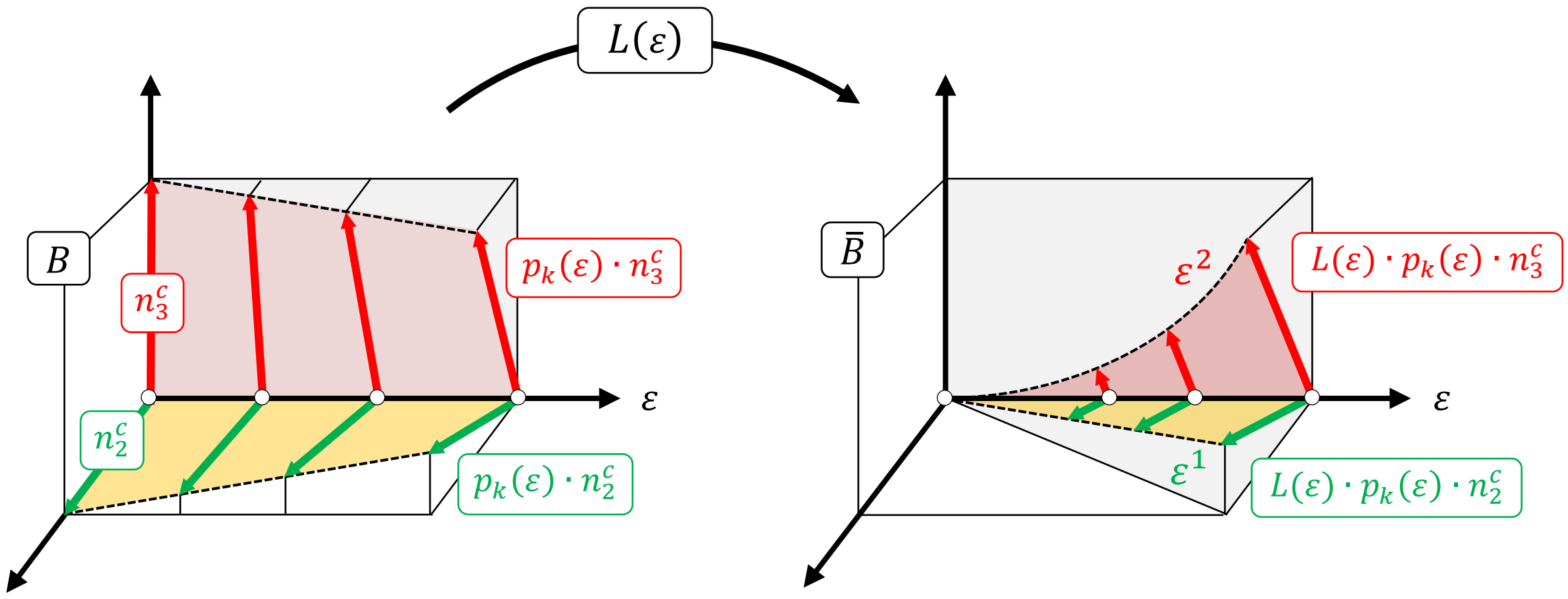

Figure 5: Blow up behaviour of $N_2^c$ and $N_3^c$ under the action of $L(\varepsilon)$.

On the left side, the domain space $B$ is indicated with elements $n_2^c \in N_2^c$ (green) and $n_3^c \in N_3^c$ (red) at $\varepsilon = 0$. Then, during variation to $\varepsilon \neq 0$, these elements are slightly perturbed by the near identity transformation $p_k(\varepsilon)$ to obtain $p_k(\varepsilon) \cdot n_2^c$ (orange surface) and $p_k(\varepsilon) \cdot n_3^c$ (red surface) respectively. Now, upon mapping of these subspaces by $L(\varepsilon)$ to the image space $\bar{B}$ on the right hand side, the image collapses to the origin at $\varepsilon = 0$, while blowing up for $\varepsilon \neq 0$ by order of $\varepsilon^1$ (orange surface) and $\varepsilon^2$ (red surface) for $n_2^c$ and $n_3^c$ respectively.

In addition, by $k$-surjectivity of the family $L(\varepsilon)$, we require the image space $\bar{B}$ to be completely build up for $\varepsilon \neq 0$ by subspaces intersecting each other transversally in the origin, as indicated in the right diagram.

Finally, the overall speed of expansion $\varepsilon^\chi$ in finite dimensions given by (2.1), can now be replaced by individual expansion rates of each of the subspaces. Due to this refinement, some of the results in finite dimensions can be improved, and above all, the technique is working with infinite dimensions as well.

For interpretation of the results in finite dimensions with respect to the generalized algebraic multiplicity $\chi$ in (2.1), application to degree theory for $\mathbb{K} = \mathbb{R}$, as well as calculation of the Milnor number $\mu$ of $ADE$-singularities using Newton-polygons, we refer to [16], [17].

Note also that in case of $L(\varepsilon) = A - \varepsilon \cdot I,\ A \in \mathbb{K}^{n,n},\ I = diag(1 \cdots 1) \in \mathbb{K}^{n,n}$, the number $\chi$ denotes the algebraic multiplicity of the eigenvalue $\varepsilon = 0$ of the matrix $A$. Aspects related to the algebraic multiplicity are treated in detail in [11] with respect to Fredholm operators and in [13] considering generalizations to infinite dimensions.

### *3. Proof of the Linearization Theorem*

First, assume surjectivity of $L_0 = G'[0]$. Then $k$-surjectivity occurs with $k = 0$ and Theorem 1 is true by (1.3)-(1.5) and $B = N_1^c \oplus N_1$, $\bar{B} = R_1$, as well as $p_0(\varepsilon) = I_B$, $S_1 = L_0$ and replacing $(\varphi_{k+1}, n_{k+1})$ by $(\varphi_1, n_1)$.

Next, proving Theorem 1 with $k \geq 1$, the Ansatz $z = z(\varepsilon) + \varepsilon^\chi \cdot r(\varepsilon) \cdot \bar{z}$ from (2.2) is replaced by

$$z = z(\varepsilon) + \varepsilon^k \cdot p_k(\varepsilon) \cdot \{ \varepsilon^k \cdot n_1^c + \cdots + \varepsilon \cdot n_k^c + (n_{k+1}^c + n_{k+1}) \} \tag{3.1}$$

$$= z(\varepsilon) + \varepsilon^k \cdot ( I_B + \varepsilon \cdot \phi_1 + \cdots + \varepsilon^k \cdot \phi_k ) \cdot \{ \varepsilon^k \cdot n_1^c + \cdots + \varepsilon \cdot n_k^c + (n_{k+1}^c + n_{k+1}) \}$$

$$= z(\varepsilon) + \varepsilon^k \cdot \{ I_B \cdot (n_{k+1}^c + n_{k+1}) + \cdots + \varepsilon^{2k} \cdot \phi_k \cdot n_1^c \} .$$

Now, it is well known from [5], [11] and [16] that in finite dimensions, we have $\chi = 1 \cdot dim\, R_2 + \cdots + k \cdot dim\, R_{k+1} \geq k$ with $R_{k+1} \neq \{0\}$, yielding $\chi = k$ merely in case of $R_2 = \cdots = R_k = \{0\}$ and $dim\, R_{k+1} = 1$. Hence, scaling by $\varepsilon^k$ instead of $\varepsilon^\chi$, typically means an increase of the cone size usable for linearization.

In addition, the variable $\bar{z}$ in $z = z(\varepsilon) + \varepsilon^\chi \cdot r(\varepsilon) \cdot \bar{z}$ is now split up into the subspaces $N_1^c, \cdots, N_k^c$ and $N_{k+1}^c, N_{k+1}$, with subspaces $N_1^c, \cdots, N_k^c$ additionally scaled from $\varepsilon^k \cdot n_1^c$ down to $\varepsilon^1 \cdot n_k^c$ for reaching an identical expansion rate of $\varepsilon^{2k}$, when mapped to the target space $\bar{B}$.

More precisely, by (2.6) and figure 5, we have $L(\varepsilon) \cdot p_k(\varepsilon) \cdot n_{i+1}^c = O(\varepsilon^i)$ for $i = 0, \ldots, k$ and scaling of $N_{i+1}^c$ by $\varepsilon^{k-i}$ implies blowing up in $\bar{B}$ by order of $\varepsilon^k = \varepsilon^{k-i} \cdot \varepsilon^i$, which finally induces blowing up by order of $\varepsilon^{2k}$, under consideration of the factor $\varepsilon^k$ in front of $p_k(\varepsilon)$ within (3.1).

In figure 6, the scaling procedure is depicted in $\mathbb{R}^n$.

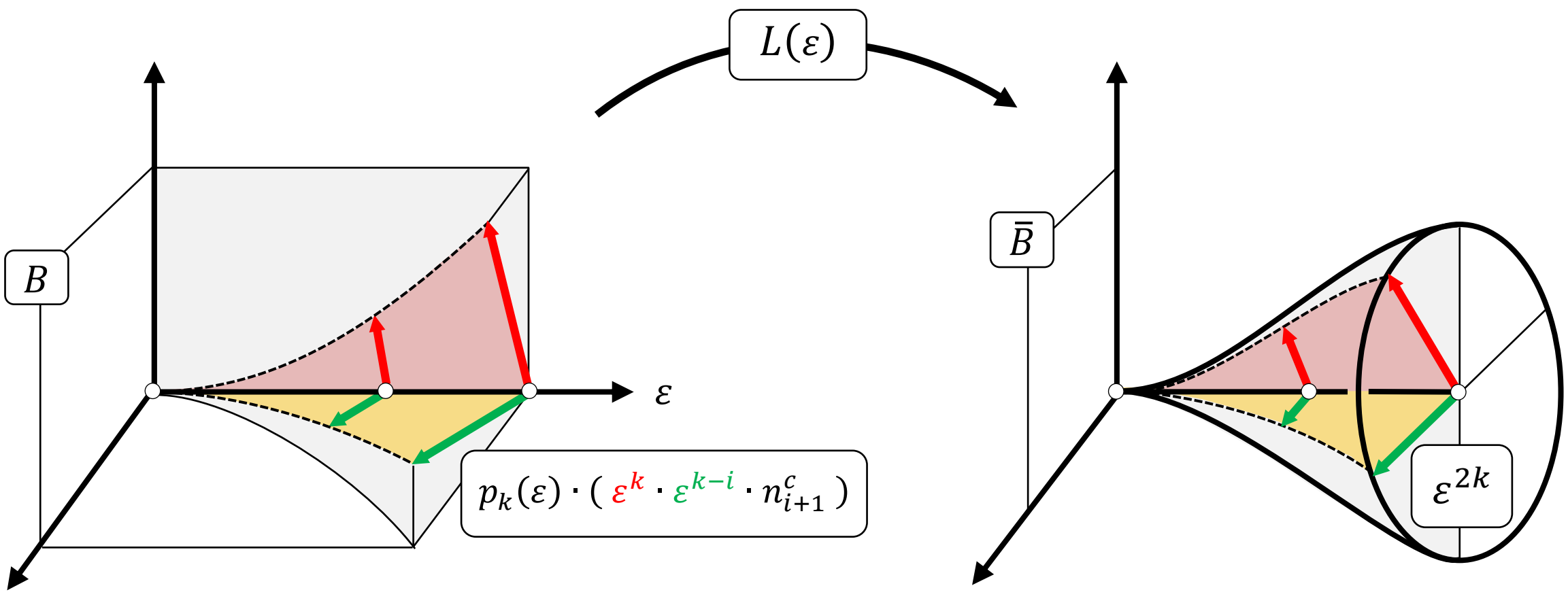

Figure 6: Uniform blow up behaviour in $\bar{B}$ of order $\varepsilon^{2k}$ after scaling of $N_1^c, \cdots, N_{k+1}^c$.

In the left diagram, the scaling of $N_{i+1}^c$ in $B$ is indicated after multiplication of $N_{i+1}^c$ with $\varepsilon^k$ and $\varepsilon^{k-i}$ (orange surface). Then, this scaled subspace is mapped by $L(\varepsilon)$ to $\bar{B}$ for $i = 0, \ldots, k$, yielding uniform blow up behaviour of order $\varepsilon^{2k}$ in the target space $\bar{B}$, i.e. a circle type blow up with radius $\varepsilon^{2k}$ occurs.

Further, note that by (3.1) and (2.6), the center line $z(\varepsilon)$ of the generalized cone is enclosed by Jordan chains with increasing rank and $\varepsilon$-starting position given by $\varepsilon^k$. In some more detail, ordering of (3.1) by powers of $\varepsilon$ implies for $k \geq 1$

$$\begin{aligned} z \;&=\; \varepsilon \cdot z_1 + \cdots + \varepsilon^{k-1} \cdot z_{k-1} \;+\; \varepsilon^k \cdot z_k \;+\; \varepsilon^{k+1} \cdot z_{k+1} \;+\; \cdots \\ &\quad + \underbrace{\varepsilon^k \cdot (n_{k+1}^c + n_{k+1})}_{start\ \ of\ Jordan\ chains} \;+\; \varepsilon^{k+1} \cdot (I_B \quad \phi_1) \cdot \begin{pmatrix} n_k^c \\ n_{k+1}^c + n_{k+1} \end{pmatrix} \;+\; \cdots \end{aligned} \tag{3.2}$$

with Jordan chains occupying powers from $\varepsilon^k$ up to $\varepsilon^{3k}$, as exemplarily and schematically indicated below in case of $k = 4$.

$$\begin{array}{cccccccccccccccc} \varepsilon & \cdots & \cdots & \varepsilon^4 & \cdots & \cdots & \cdots & \varepsilon^8 & \cdots & \cdots & \cdots & \varepsilon^{12} & \cdots & \cdots & \cdots & \cdots \\ z_1 & z_2 & z_3 & z_4 & z_5 & z_6 & z_7 & z_8 & z_9 & z_{10} & z_{11} & z_{12} & z_{13} & z_{14} & z_{15} & \cdots \\ & & & \boxtimes & \boxtimes & \boxtimes & \boxtimes & \boxtimes & \boxtimes & \boxtimes & \boxtimes & \boxtimes & & & & \\ & left \leftarrow & & & & & Jordan\ chains & & & & & \rightarrow right & & & & \end{array} \tag{3.3}$$

Note also that in (3.1), two polynomials of degree $k$ in $\varepsilon$ are multiplied by $p_k(\varepsilon) \cdot \{\varepsilon^k \cdot n_1^c + \cdots + n_{k+1}\}$, yielding a polynomial of degree $\varepsilon^{2k}$, which is finally multiplied in (3.1) with $\varepsilon^k$, determining the starting position of the chains within the Ansatz. In (3.3), the positions of the chains from $\varepsilon^k$ up to $\varepsilon^{2k}$ are labeled by red colour, which are essential for the procedure to work, due to maximal length $k$ of Jordan chains causing $k$-surjectivity. In [17], only this minimal constellation is used and black marked positions between $\varepsilon^{2k+1}$ and $\varepsilon^{3k}$ are set to zero.

The starting position of the Jordan chains may be shifted to the left or to the right, simply by replacing red marked $\varepsilon^k$ in (3.1) by $\varepsilon^{k \pm l}$, $l \geq 1$. By a shift to the right according to $\varepsilon^{k+l}$, we obtain uniform expansion in $\bar{B}$ by $\varepsilon^{2k+l}$, allowing linearization and calculation of solutions with $\varepsilon^{2k}$ replaced by $\varepsilon^{2k+l}$ in Theorem 1 and Corollaries 1, 2. In other words, a shift to the right merely narrows the cone and effectively no improvement is achieved.

On contrary, by a shift to the left according to $\varepsilon^{k-l}$, a complete sequence of Newton Lemma type results are obtained in a rather systematic way, based on appropriate assumptions with respect to higher order terms in (1.1). This investigation is performed in section 6.

Now, plugging (3.1) into the expansion (1.1), we see from (2.6), (2.7)

$$\begin{aligned} G[*] \coloneqq\; & G[\, z(\varepsilon) \;+\; \overbrace{\varepsilon^k \cdot p_k(\varepsilon) \cdot \{\, \varepsilon^k \cdot n_1^c \;+ \cdots +\; \varepsilon \cdot n_k^c \;+\; n_{k+1}^c + n_{k+1} \,\}}^{=:\, \varepsilon^k \cdot b} \,] \\ =\; & G[\, z(\varepsilon) \,] \;+\; L(\varepsilon) \cdot \varepsilon^k \cdot p_k(\varepsilon) \cdot \{\, \varepsilon^k \cdot n_1^c \;+ \cdots +\; \varepsilon \cdot n_k^c \;+\; n_{k+1}^c + n_{k+1} \,\} \\ & + \mathcal{B}(\varepsilon) \cdot \left(\varepsilon^k \cdot b\right)^2 \;+\; r\left(\varepsilon, \varepsilon^k \cdot b\right) \cdot \left(\varepsilon^k \cdot b\right)^3 \end{aligned} \tag{3.4}$$

$$= G[\,z(\varepsilon)\,] \;+\; \varepsilon^k \cdot \{\, \overbrace{\varepsilon^k \cdot L(\varepsilon)\cdot p_k(\varepsilon)\cdot n_1^c}^{=\,S_1 n_1^c + O(\varepsilon)} \;+\cdots+\; \overbrace{L(\varepsilon)\cdot p_k(\varepsilon)\cdot n_{k+1}^c}^{=\,\varepsilon^k\cdot S_{k+1} n_{k+1}^c + O(\varepsilon^{k+1})} \;+\; \overbrace{L(\varepsilon)\cdot p_k(\varepsilon)\cdot n_{k+1}}^{=\,O(\varepsilon^{k+1})} \,\}$$

$$+\; \varepsilon^{2k}\cdot \mathcal{B}(\varepsilon)\cdot b^2 \;+\; \varepsilon^{3k}\cdot r\big(\varepsilon, \varepsilon^k\cdot b\big)\cdot b^3$$

$$= G[\,z(\varepsilon)\,] \;+\; \varepsilon^{2k}\cdot\{\, S_1\cdot n_1^c \;+\cdots+\; S_{k+1}\cdot n_{k+1}^c + O(\varepsilon)\,\}$$

$$+\; \varepsilon^{2k}\cdot[\,\mathcal{B}_0 + O(\varepsilon)\,]\cdot b^2 \;+\; \varepsilon^{3k}\cdot r\big(\varepsilon, \varepsilon^k\cdot b\big)\cdot b^3$$

and collecting powers of $\varepsilon^{2k}$, we end up with

$$G[*] = G[\,z(\varepsilon)\,] \;+\; \varepsilon^{2k}\cdot\{\,(S_1\cdots S_{k+1})\cdot \overbrace{\begin{pmatrix} n_1^c \\ \vdots \\ n_{k+1}^c \end{pmatrix}}^{=\,n^c} \tag{3.5}$$

$$+\; \mathcal{B}_0\cdot(\,n_{k+1}^c + n_{k+1}\,)^2 \;+\; \varepsilon\cdot R(\,\varepsilon, n^c, n_{k+1}\,)\,\}$$

in case of $k \geq 1$. Here, $R(\cdot)$ represents a smooth family of remainder operators satisfying $R(\varepsilon, 0{,}0) = 0$. In addition, we see that the operator $\mathcal{B}_0 = \frac{1}{2} G''[0]$ is also part of the analysis. In section 6 of extended Newton Lemmas, further derivatives of $G$ at the origin will be taken into consideration.

Now, for turning (3.5) into the linearization (1.14) of Theorem 1 (i), we merely have to eliminate the summands $\mathcal{B}_0\cdot(n_{k+1}^c + n_{k+1})^2 + \varepsilon\cdot R(\varepsilon, n^c, n_{k+1})$ in (3.5) by solving a smooth remainder equation of the form $H(\varepsilon, n^c, n_{k+1}, \varphi) = 0$, $H: U \times N^c \times N_{k+1} \times N^c \to \bar{B}$ defined by

$$(S_1\cdots S_{k+1})\cdot n^c \;+\; \mathcal{B}_0\cdot(\,n_{k+1}^c + n_{k+1}\,)^2 \;+\; \varepsilon\cdot R(\varepsilon, n^c, n_{k+1}) \;=\; (S_1\cdots S_{k+1})\cdot\varphi$$

$$\Leftrightarrow \quad H(\varepsilon, n^c, n_{k+1}, \varphi) \tag{3.6}$$

$$:= (S_1\cdots S_{k+1})\cdot(n^c - \varphi) \;+\; \mathcal{B}_0\cdot(n_{k+1}^c + n_{k+1})^2 \;+\; \varepsilon\cdot R(\varepsilon, n^c, n_{k+1}) = 0\,.$$

We aim to solve $H(\varepsilon, n^c, n_{k+1}, \varphi) = 0$ with respect to $n^c \in N^c = N_1^c \times \cdots \times N_{k+1}^c$. Obviously,

$$H(\varepsilon, 0{,}0{,}0) = 0 \quad \textit{and} \quad H_{n^c}(0{,}0{,}0{,}0) = (S_1\cdots S_{k+1}) \in GL[\,N^c, \bar{B}\,]\,, \tag{3.7}$$

and hence by standard implicit function theorem, the existence of a locally unique and smooth operator function $\psi^c(\varepsilon, \varphi, n_{k+1})$, satisfying

$$H[\,\varepsilon,\; \psi^c(\varepsilon, \varphi, n_{k+1}),\; n_{k+1},\; \varphi\,] = 0 \tag{3.8}$$

and

$$G[\,z(\varepsilon) \;+\; \varepsilon^k\cdot p_k(\varepsilon)\cdot\{\, \overbrace{\varepsilon^k\cdot\psi_1^c(\varepsilon,\varphi,n_{k+1}) \;+\cdots+\; \psi_{k+1}^c(\varepsilon,\varphi,n_{k+1})}^{=\,\psi^c(\,\varepsilon,\varphi,n_{k+1}\,)} \;+\; n_{k+1}\,\}\,] \tag{3.9}$$

$$\underset{(3.5)}{\overset{(3.4)}{=}} G[\,z(\varepsilon)\,] \;+\; \varepsilon^{2k}\cdot(S_1\cdots S_{k+1})\cdot\varphi\,,$$

is assured, implying local linearization (1.14), (1.15) as desired. Tringularity of $P_k(\varepsilon)$ follows from (2.8). Next, concerning the properties of $\psi^c(\varepsilon,\varphi,n_{k+1})$ at $\varepsilon = 0$ in (1.16), consider splitting of the equation by

$$H(\varepsilon,n^c,n_{k+1},\varphi) = 0 \in \bar{B} \quad \Leftrightarrow \quad \mathcal{P}_i \cdot H(\varepsilon,n^c,n_{k+1},\varphi) = 0 \in R_i\,, \tag{3.10}$$

$i = 1,\cdots,k+1$, based on $\bar{B} = R_1 \oplus \cdots \oplus R_{k+1}$ from (1.10). Then, due to $S_i \in GL[N_i^c, R_i\,]$, equations (3.6), (3.10) simplify according to

$$\begin{pmatrix} S_1 & & \\ & \ddots & \\ & & S_k \end{pmatrix} \cdot \begin{pmatrix} n_1^c - \varphi_1 \\ \vdots \\ n_k^c - \varphi_k \end{pmatrix} + \begin{pmatrix} \mathcal{P}_1 \\ \vdots \\ \mathcal{P}_k \end{pmatrix} \cdot [\,\mathcal{B}_0 \cdot (\,n_{k+1}^c + n_{k+1}\,)^2 \;+\; \varepsilon \cdot R(\varepsilon,n^c,n_{k+1})\,] = 0 \tag{3.11}$$

$$S_{k+1} \cdot (n_{k+1}^c - \varphi_{k+1}) + \mathcal{P}_{k+1} \cdot [\,\mathcal{B}_0 \cdot (\,n_{k+1}^c + n_{k+1}\,)^2 \;+\; \varepsilon \cdot R(\varepsilon,n^c,n_{k+1})\,] = 0\,, \tag{3.12}$$

implying decoupling of the equations at $\varepsilon = 0$ in the sense that the last equation (3.12) is independent of $n_1^c,\ldots,n_k^c$. Hence, we can solve (3.12) at $\varepsilon = 0$ with respect to $n_{k+1}^c$, then plugging these solutions into equations (3.11) for calculation of $n_1^c,\ldots,n_k^c$ at $\varepsilon = 0$ from pure linear relations. First, define

$$Q(n_{k+1}^c,n_{k+1},\varphi_{k+1}) := S_{k+1} \cdot (n_{k+1}^c - \varphi_{k+1}) + \mathcal{P}_{k+1} \cdot \mathcal{B}_0 \cdot (\,n_{k+1}^c + n_{k+1}\,)^2\,, \tag{3.13}$$

satisfying $Q(0,0,0) = 0$ and

$$Q_{n_{k+1}^c}(n_{k+1}^c,n_{k+1},\varphi_{k+1}) = S_{k+1} + 2 \cdot \mathcal{P}_{k+1} \cdot \mathcal{B}_0 \cdot (\,n_{k+1}^c + n_{k+1}\,) \in GL[\,N_{k+1}^c, R_{k+1}\,] \tag{3.14}$$

for $n_{k+1}^c, n_{k+1}$ chosen sufficiently small, under consideration of $S_{k+1} \in GL[N_{k+1}^c, R_{k+1}]$ and continuity of $\mathcal{P}_{k+1}$ and $\mathcal{B}_0$. Hence, there exists a locally unique mapping $n_{k+1}^c(n_{k+1},\varphi_{k+1})$ satisfying

$$n_{k+1}^c(0,0) = 0 \qquad and \qquad Q[\,n_{k+1}^c(n_{k+1},\varphi_{k+1}),\; n_{k+1},\; \varphi_{k+1}\,] = 0\,. \tag{3.15}$$

Further, by Taylor's formula and implicit differentiation of (3.15), we obtain

$$\begin{aligned} n_{k+1}^c(n_{k+1},\varphi_{k+1}) \;&=\; n_{k+1}^c(0,0) \;+\; {n_{k+1}^c}'(0,0) \cdot \begin{pmatrix} n_{k+1} \\ \varphi_{k+1} \end{pmatrix} \;+\; \chi_{k+1}^c(n_{k+1},\varphi_{k+1}) \cdot \begin{pmatrix} n_{k+1} \\ \varphi_{k+1} \end{pmatrix}^2 \\ &\overset{(3.15)}{=} 0 \;+\; [\,0 \quad I_{N^c}\,] \cdot \begin{pmatrix} n_{k+1} \\ \varphi_{k+1} \end{pmatrix} \;+\; \chi_{k+1}^c(n_{k+1},\varphi_{k+1}) \cdot \begin{pmatrix} n_{k+1} \\ \varphi_{k+1} \end{pmatrix}^2 \\ &=\; \varphi_{k+1} \;+\; \underbrace{\chi_{k+1}^c(n_{k+1},\varphi_{k+1})}_{\to\, N_{k+1}^c} \cdot \begin{pmatrix} n_{k+1} \\ \varphi_{k+1} \end{pmatrix}^2 \\ &=\; \psi_{k+1}^c(0,\varphi,n_{k+1}) \end{aligned} \tag{3.16}$$

with smooth remainder function $\chi_{k+1}^c(n_{k+1},\varphi_{k+1})$ mapping to $N_{k+1}^c$. Next, plugging $n_{k+1}^c(n_{k+1},\varphi_{k+1})$ into (3.11) at $\varepsilon = 0$, we arrive at the linear system

$$\begin{pmatrix} S_1 & & \\ & \ddots & \\ & & S_k \end{pmatrix} \cdot \begin{pmatrix} n_1^c - \varphi_1 \\ \vdots \\ n_k^c - \varphi_k \end{pmatrix} + \begin{pmatrix} \mathcal{P}_1 \\ \vdots \\ \mathcal{P}_k \end{pmatrix} \cdot \mathcal{B}_0 \cdot [\, n_{k+1}^c(n_{k+1}, \varphi_{k+1}) + n_{k+1} \,]^2 \; = 0 \tag{3.17}$$

with solutions explicitly given by

$$\begin{array}{lcccl} n_1^c(\varphi_1, n_{k+1}, \varphi_{k+1}) & = & \varphi_1 & - & S_1^{-1} \cdot \mathcal{P}_1 \cdot \mathcal{B}_0 \cdot [\, n_{k+1}^c(n_{k+1}, \varphi_{k+1}) + n_{k+1} \,]^2 \\ \vdots & & \vdots & \vdots & \\ n_k^c(\varphi_k, n_{k+1}, \varphi_{k+1}) & = & \varphi_k & - & S_k^{-1} \cdot \mathcal{P}_k \cdot \mathcal{B}_0 \cdot [\, n_{k+1}^c(n_{k+1}, \varphi_{k+1}) + n_{k+1} \,]^2 \,. \end{array} \tag{3.18}$$

Now, by (3.16) and (3.18), we end up with the final form for $i = 1, \ldots, k$

$$\begin{aligned} \psi_i^c(0, \varphi, n_{k+1}) = n_i^c(\varphi_i, n_{k+1}, \varphi_{k+1}) &= \varphi_i - S_i^{-1} \cdot \mathcal{P}_i \cdot \mathcal{B}_0 \cdot [\, n_{k+1}^c(n_{k+1}, \varphi_{k+1}) + n_{k+1} \,]^2 \\ &= \varphi_i - S_i^{-1} \cdot \mathcal{P}_i \cdot \mathcal{B}_0 \cdot [\; \varphi_{k+1} \; + \; \chi_{k+1}^c(n_{k+1}, \varphi_{k+1}) \cdot \begin{pmatrix} n_{k+1} \\ \varphi_{k+1} \end{pmatrix}^2 + n_{k+1} \,]^2 \\ &=: \varphi_i \; + \; \underbrace{\chi_i^c(n_{k+1}, \varphi_{k+1})}_{\to\, N_i^c} \cdot \begin{pmatrix} n_{k+1} \\ \varphi_{k+1} \end{pmatrix}^2 \end{aligned} \tag{3.19}$$

with smooth operator functions $\chi_i^c(n_{k+1}, \varphi_{k+1})$ mapping to $N_i^c$. Finally, using $H(\varepsilon, 0,0,0) = 0$ from (3.7), uniqueness shows $\psi_1^c(\varepsilon, 0,0) = \cdots = \psi_{k+1}^c(\varepsilon, 0,0) = 0$ and Theorem 1 (i) is proved.

Concerning (ii), we have to show that for fixed $\varepsilon \neq 0$, the mapping $p_k(\varepsilon) \cdot \{\Psi^c(\varepsilon, \varphi, n_{k+1}) + n_{k+1}\}$ represents a diffeomorphism in $B$ between

$$V = \{\, (\varphi, n_{k+1}) \in N^c \times N_{k+1} \;\mid\; \|(\varphi, n_{k+1})\| \ll 1 \,\} \tag{3.20}$$

and an open neighbourhood $V_0$ of the origin $0 \in B$.

First, the near identity mapping $p_k(\varepsilon) = I_B + \varepsilon \cdot \phi_1 + \cdots + \varepsilon^k \cdot \phi_k$ represents an isomorphism in $B$ for $\varepsilon$ sufficiently small and it remains to show diffeomorphic behaviour of the nonlinear mapping $\Psi^c(\varepsilon, \varphi, n_{k+1}) + n_{k+1}$, which can obviously be written as a product according to

$$\begin{aligned} &\Psi^c(\, \varepsilon, \varphi, n_{k+1} \,) + n_{k+1} \\ &= \varepsilon^k \cdot \psi_1^c(\varepsilon, \varphi, n_{k+1}) \; + \cdots + \; \varepsilon^1 \cdot \psi_k^c(\varepsilon, \varphi, n_{k+1}) \; + \; \psi_{k+1}^c(\varepsilon, \varphi, n_{k+1}) \; + \; n_{k+1} \\ &= (\, \varepsilon^k \cdot P_1 \; + \cdots + \; \varepsilon^1 \cdot P_k \; + \; P_{k+1} \; + \; P_{N_{k+1}} \,) \cdot [\, \psi_1^c(\varepsilon, \varphi, n_{k+1}) \; + \cdots + \; \psi_{k+1}^c(\varepsilon, \varphi, n_{k+1}) \; + \; n_{k+1} \,] \end{aligned} \tag{3.21}$$

with $P_1, \cdots, P_{k+1}$ and $P_{N_{k+1}}$ denoting continous projection to $N_1^c, \cdots, N_{k+1}^c$ and $N_{k+1}$ respectively. Now, the first factor represents a scaled version of the identity map in $B$, as well as an isomorphism of $B$ for $\varepsilon \neq 0$.

Concerning the nonlinear second factor, we calculate the inverse function simply by solving the equation

$$\psi_1^c(\varepsilon, \varphi, n_{k+1}) \; + \cdots + \; \psi_{k+1}^c(\varepsilon, \varphi, n_{k+1}) \; + \; n_{k+1} = b$$

$$\Leftrightarrow \quad D(\varepsilon,\varphi,n_{k+1},b) := \psi_1^c(\varepsilon,\varphi,n_{k+1}) + \cdots + \psi_{k+1}^c(\varepsilon,\varphi,n_{k+1}) + n_{k+1} - b = 0\,. \tag{3.22}$$

We have $D(0,0,0,0) = 0$ and we see from (3.16) and (3.19)

$$D(0,\varphi,n_{k+1},b) = \psi_1^c(0,\varphi,n_{k+1}) + \cdots + \psi_{k+1}^c(0,\varphi,n_{k+1}) + n_{k+1} - b \tag{3.23}$$

$$= \varphi_1 + \chi_1^c(n_{k+1},\varphi_{k+1}) \cdot \begin{pmatrix} n_{k+1} \\ \varphi_{k+1} \end{pmatrix}^2 + \cdots + \varphi_{k+1} + \chi_{k+1}^c(n_{k+1},\varphi_{k+1}) \cdot \begin{pmatrix} n_{k+1} \\ \varphi_{k+1} \end{pmatrix}^2 + n_{k+1} - b$$

$$=: (\varphi_1 + \cdots + \varphi_{k+1} + n_{k+1}) + \chi^c(n_{k+1},\varphi_{k+1}) \cdot \begin{pmatrix} n_{k+1} \\ \varphi_{k+1} \end{pmatrix}^2 - b\,,$$

implying regularity in the origin according to

$$D_{(\varphi,n_{k+1})}(0,0,0,0) = P_1 + \cdots + P_{k+1} + P_{N_{k+1}} = I_B \in GL[B,B]\,. \tag{3.24}$$

Hence, there exist open neighbourhoods $U \subset \mathbb{K}$, $V \subset N^c \times N_{k+1}$, $V_0 \subset B$ of the origins and a smooth unique mapping $(\varphi,n_{k+1})(\varepsilon,b)$ such that every solution of $D(\varepsilon,\varphi,n_{k+1},b) = 0$ in $U \times V \times V_0$ is given by $(\varphi,n_{k+1})(\varepsilon,b)$. In particular, for $\varepsilon$ fixed, the inverse mapping of

$$\psi_1^c(\varepsilon,\varphi,n_{k+1}) + \cdots + \psi_{k+1}^c(\varepsilon,\varphi,n_{k+1}) + n_{k+1} \tag{3.25}$$

is given by $(\varphi,n_{k+1})(\varepsilon,b)$ for $b \in V_0$, yielding diffeomorphic behaviour between $V$ and $V_0$ and Theorem 1 (ii) is shown.

## *4. Proof of Corollary 1*

First, for $k \geq 0$ we obtain from (1.24), (1.25) and Theorem 1 (i)

$$G[\,z(\varepsilon) + \varepsilon^k \cdot \{\Psi^c(\varepsilon,\varphi + \varphi_0(\varepsilon),n_{k+1}) + n_{k+1}\}\,]$$

$$= G[\,z(\varepsilon)\,] + \varepsilon^{2k} \cdot (S_1 \cdots S_{k+1}) \cdot [\,\varphi + \varphi_0(\varepsilon)\,]$$

$$= \varepsilon^{2k} \cdot [\,\bar{b}(\varepsilon) + (S_1 \cdots S_{k+1}) \cdot \varphi + \underbrace{(S_1 \cdots S_{k+1}) \cdot \varphi_0(\varepsilon)}_{= -\bar{b}(\varepsilon)}\,] \tag{4.1}$$

$$= \varepsilon^{2k} \cdot (S_1 \cdots S_{k+1}) \cdot \varphi\,.$$

Further, smallness of $\varphi_0(\varepsilon) = -(S_1 \cdots S_{k+1})^{-1} \cdot \bar{b}(\varepsilon) \in N^c$ is assured by the second condition $\|\bar{b}(0)\| \ll 1$ in (1.24) and Corollary 1 (i) is shown.

Concerning (ii), assume $k \geq 1$. Again, from Theorem 1 (i), we obtain

$$G[\,z(\varepsilon) + \varepsilon^k \cdot \{\Psi^c(\varepsilon,\varphi,n_{k+1}) + n_{k+1}\}\,]$$

$$= G[\,z(\varepsilon)\,] + \varepsilon^{2k} \cdot (S_1 \cdots S_{k+1}) \cdot \varphi \tag{4.2}$$

$$= \varepsilon^q \cdot [\,\underbrace{\bar{b}(\varepsilon)}_{\neq 0} + \varepsilon^{2k-q} \cdot (S_1 \cdots S_{k+1}) \cdot \varphi\,] \neq 0$$

for $\varepsilon \in U \setminus \{0\}$, $(\varphi,n_{k+1}) \in V$ and $U$ chosen sufficiently small.

## 5. Proof of Newton's Lemma

The proof of Corollary 2 starts similarly to the proof of Theorem 1. From (3.4), (3.5) and assumed approximation of order $2k$ by $G[z(\varepsilon)] = \varepsilon^{2k} \cdot \bar{b}(\varepsilon)$, we obtain for $k \geq 1$

$$G[*] = G[\, z(\varepsilon) + \varepsilon^k \cdot p_k(\varepsilon) \cdot \{ \, \varepsilon^k \cdot n_1^c \, + \cdots + \, \varepsilon \cdot n_k^c \, + \, n_{k+1}^c + n_{k+1} \, \} \,] \tag{5.1}$$

$$\overset{(3.5)}{=} G[\, z(\varepsilon)\,] \, + \, \varepsilon^{2k} \cdot \{ \, (S_1 \cdots S_{k+1}) \cdot n^c + \mathcal{B}_0 \cdot (\, n_{k+1}^c + n_{k+1}\,)^2 \, + \, \varepsilon \cdot R(\, \varepsilon, n^c, n_{k+1}\,) \, \}$$

$$= \varepsilon^{2k} \cdot \, \{ \, \bar{b}(\varepsilon) \, + \, (\, S_1 \cdots S_{k+1}\,) \cdot n^c \, + \, \mathcal{B}_0 \cdot (\, n_{k+1}^c + n_{k+1}\,)^2 \, + \, \varepsilon \cdot R(\, \varepsilon, n^c, n_{k+1}\,) \, \}$$

and the blown up remainder equation to be solved reads

$$(\, S_1 \cdots S_{k+1}\,) \cdot n^c \, + \, \mathcal{B}_0 \cdot (\, n_{k+1}^c + n_{k+1}\,)^2 \, + \, \bar{b}(\varepsilon) \, + \, \varepsilon \cdot R(\, \varepsilon, n^c, n_{k+1}\,) = 0\,. \tag{5.2}$$

As previously, the equation is first solved at $\varepsilon = 0$ and then continued to $\varepsilon \neq 0$. Splitting, based on $\bar{B} = R_1 \oplus \cdots \oplus R_{k+1}$, implies at $\varepsilon = 0$

$$\begin{array}{lllllll} S_1 \cdot n_1^c & + & \mathcal{P}_1 \cdot \mathcal{B}_0 \cdot (\, n_{k+1}^c + n_{k+1}\,)^2 & + & \mathcal{P}_1 \cdot \bar{b}(0) & = & 0 \\ & \vdots & & \vdots & & & \\ S_k \cdot n_k^c & + & \mathcal{P}_k \cdot \mathcal{B}_0 \cdot (\, n_{k+1}^c + n_{k+1}\,)^2 & + & \mathcal{P}_k \cdot \bar{b}(0) & = & 0 \\ S_{k+1} \cdot n_{k+1}^c & + & \mathcal{P}_{k+1} \cdot \mathcal{B}_0 \cdot (\, n_{k+1}^c + n_{k+1}\,)^2 & + & \mathcal{P}_{k+1} \cdot \bar{b}(0) & = & 0 \end{array} \tag{5.3}$$

with quadratic last equation concerning $n_{k+1}^c$ and linear equations with respect to $n_1^c, \cdots, n_k^c$. Our aim consists in stating conditions, which ensure the existence of a regular solution of the last equation. Different conditions will distinguish the cases (i) and (ii) of Corollary 2.

First, assume $\|\mathcal{P}_{k+1} \cdot \bar{b}(0)\| \ll 1$, as stated in Corollary 2 (i). We examine the equation

$$H(n_{k+1}^c, n_{k+1}, r_{k+1}) \coloneqq S_{k+1} \cdot n_{k+1}^c + \mathcal{P}_{k+1} \cdot \mathcal{B}_0 \cdot (\, n_{k+1}^c + n_{k+1}\,)^2 + r_{k+1} = 0 \tag{5.4}$$

with smooth mapping $H\colon N_{k+1}^c \times N_{k+1} \times R_{k+1} \to R_{k+1}$, satisfying

$$H(0,0,0) = 0 \quad and \quad H_{n_{k+1}^c}(0,0,0) = S_{k+1} \in GL[\, N_{k+1}^c, R_{k+1}\,]\,. \tag{5.5}$$

Hence, in a neighbourhood of the origin, there exists a unique mapping $n_{k+1}^c(n_{k+1}, r_{k+1})$ with $H[n_{k+1}^c(n_{k+1}, r_{k+1}), n_{k+1}, r_{k+1}] = 0$ and by setting $\bar{n}_{k+1}^c(0, n_{k+1}) \coloneqq n_{k+1}^c(n_{k+1}, \mathcal{P}_{k+1} \cdot \bar{b}(0))$, we end up with the last equation of (5.3) solved by

$$S_{k+1} \cdot \bar{n}_{k+1}^c(0, n_{k+1}) + \mathcal{P}_{k+1} \cdot \mathcal{B}_0 \cdot [\, \bar{n}_{k+1}^c(0, n_{k+1}) + n_{k+1}\,]^{\,2} + \underbrace{\mathcal{P}_{k+1} \cdot \bar{b}(0)}_{\|\cdot\| \ll 1} = 0\,. \tag{5.6}$$

Further, from (5.3), the formulas concerning first $k$ components $\bar{n}_i^c(0, n_{k+1})\,, i = 1, \cdots, k$ follow immediately, as stated in (1.33). Now, the linearization of the blown up remainder equation (5.2) with respect to $n^c$ at these solutions satisfies

$$(\, S_1 \cdots S_k \;|\; S_{k+1} + 2 \cdot \mathcal{B}_0 \cdot [\, \bar{n}_{k+1}^c(0, n_{k+1}) + n_{k+1}\,]\,) \in GL[\, N^c, B\,] \tag{5.7}$$

for $\|n_{k+1}\| \ll 1$ and continuation of $\bar{n}_i^c(0, n_{k+1})$ to $\bar{n}_i^c(\varepsilon, n_{k+1})$ is ensured by implicit function theorem for $\varepsilon \in U$ and $i = 1, \cdots, k+1$.

Next, concerning Corollary 2 (ii), the second summand in (5.2) is provided with a suitable assumption for guaranteeing regular solutions of the remainder equation. Here, the requirement $\mathcal{B}_0 \cdot (n_{k+1}^c)^2 = 0$ is turning the last equation in (5.3) into a linear equation with respect to $n_{k+1}^c$ according to

$$S_{k+1} \cdot n_{k+1}^c \; + \; \mathcal{P}_{k+1} \cdot \mathcal{B}_0 \cdot [\, (n_{k+1}^c)^2 + 2 \cdot n_{k+1}^c \cdot n_{k+1} + (n_{k+1})^2 \,] \; + \; \mathcal{P}_{k+1} \cdot \bar{b}(0)$$

$$= (\, S_{k+1} + 2 \cdot \mathcal{P}_{k+1} \cdot \mathcal{B}_0 \cdot n_{k+1} \,) \cdot n_{k+1}^c \; + \; \mathcal{P}_{k+1} \cdot \mathcal{B}_0 \cdot (n_{k+1})^2 \; + \; \mathcal{P}_{k+1} \cdot \bar{b}(0) = 0 \qquad (5.8)$$

with corresponding solution manifold $\bar{n}_{k+1}^c(0, n_{k+1})$ given by (1.36). Note that $(S_{k+1} + 2 \cdot \mathcal{P}_{k+1} \cdot \mathcal{B}_0 \cdot n_{k+1})$ belongs to $GL[N_{k+1}^c, R_{k+1}]$ for $\|n_{k+1}\| \ll 1$. First $k$ components $\bar{n}_i^c(0, n_{k+1}), i = 1, \cdots, k$ follow again from (5.3).

Finally, under consideration of $\mathcal{B}_0 \cdot (n_{k+1}^c)^2 = 0$, the linearization with respect to $n^c$ of the remainder equation (5.2) satisfies at $\varepsilon = 0$

$$(\, S_1 \cdots S_k \;|\; S_{k+1} + 2 \cdot \mathcal{B}_0 \cdot n_{k+1} \,) \in GL[\, N^c, B \,] \qquad (5.9)$$

for $\|n_{k+1}\| \ll 1$ and continuation of $\bar{n}_i^c(0, n_{k+1})$ to $\bar{n}_i^c(\varepsilon, n_{k+1})$ is again ensured by implicit function theorem for $\varepsilon \in U$ and $i = 1, \cdots, k+1$, thus finishing the proof of Corollary 2.

## *6. Improving Newton's Lemma by shift of Jordan chains*

The basic idea of the paper consists in enclosing the center line $z(\varepsilon)$ by Jordan chains of increasing rank obtained from the family of linear mappings $L(\varepsilon) = G'[z(\varepsilon)]$, which is supposed to be $k$-surjective. Up to now, the starting position of the Jordan chains within the Ansatz was given by $\varepsilon^k$ according to (3.1)-(3.3). In this section, the starting position will be shifted to the left, thereby creating in a systematic way conditions to be satisfied for the new Ansatz to work.

Note that shifting Jordan chains to the left, means increase of the cone and the results of the linearization Theorem 1 and Corollaries 1, 2 will be valid within an enlarged region of $B$, as indicated in the figure below within $\mathbb{R}^n$.

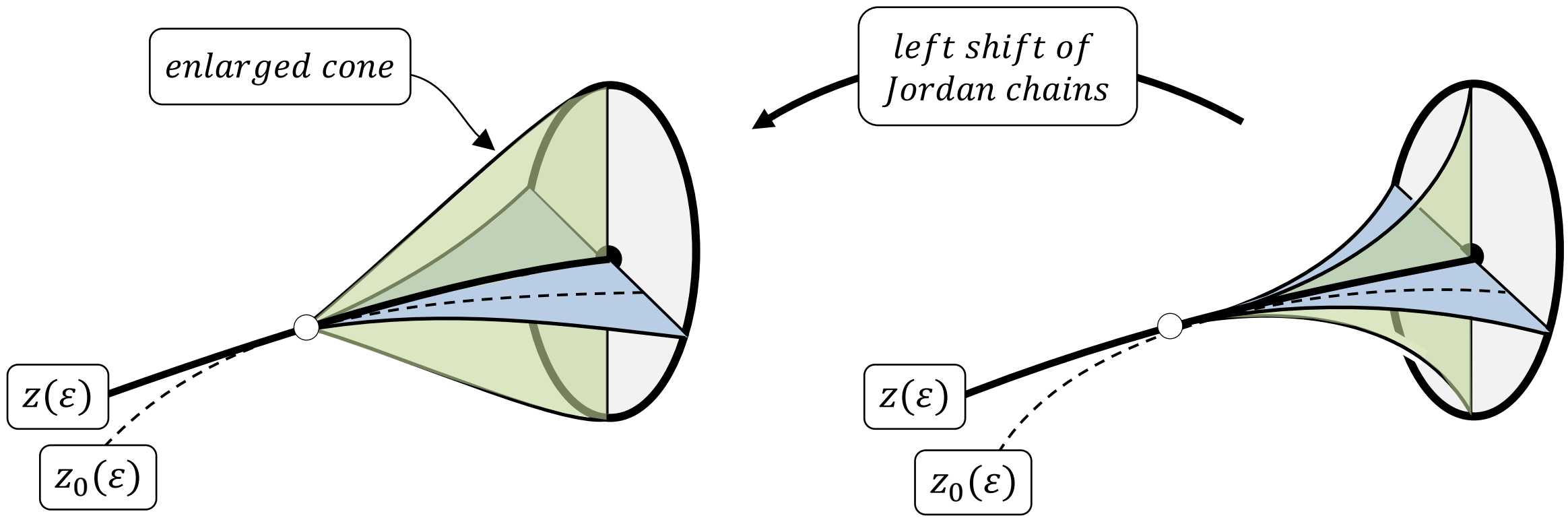

Figure 7: Increase of the cone by left shift of Jordan chains.

We also indicate possible solutions $z_0(\varepsilon)$ of $G[z] = 0$ by dashed lines, occuring within the cones, if appropriate approximation requirements concerning the center line $z(\varepsilon)$ are additionally satisfied. For simplicity, only one half cone is shown in each diagram.

First, for $k \geq 2$, the cone is shifted to the left merely by one position according to

$$z = z(\varepsilon) + \underbrace{\varepsilon^{k-1}}_{\substack{starting\\position}} \cdot p_k(\varepsilon) \cdot \{ \varepsilon^k \cdot n_1^c + \cdots + \varepsilon \cdot n_k^c + n_{k+1}^c + n_{k+1} \}, \tag{6.1}$$

implying from expansion (1.1)

$$\begin{aligned} G[*] &:= G[\, z(\varepsilon) + \overbrace{\varepsilon^{k-1} \cdot p_k(\varepsilon) \cdot \{ \varepsilon^k \cdot n_1^c + \cdots + \varepsilon \cdot n_k^c + n_{k+1}^c + n_{k+1} \}}^{=:\, \varepsilon^{k-1} \cdot b} \,] \\ &= G[\, z(\varepsilon) \,] + L(\varepsilon) \cdot \varepsilon^{k-1} \cdot p_k(\varepsilon) \cdot \{ \varepsilon^k \cdot n_1^c + \cdots + \varepsilon \cdot n_k^c + n_{k+1}^c + n_{k+1} \} \\ &\quad + \mathcal{B}(\varepsilon) \cdot \left(\varepsilon^{k-1} \cdot b\right)^2 + r\left(\varepsilon, \varepsilon^{k-1} \cdot b\right) \cdot \left(\varepsilon^{k-1} \cdot b\right)^3 \end{aligned} \tag{6.2}$$

$$\begin{aligned} &= G[\, z(\varepsilon) \,] + \varepsilon^{k-1} \cdot \{ \varepsilon^k \cdot \overbrace{L(\varepsilon) \cdot p_k(\varepsilon) \cdot n_1^c}^{= S_1 n_1^c + O(\varepsilon)} + \cdots + \overbrace{L(\varepsilon) \cdot p_k(\varepsilon) \cdot n_{k+1}^c}^{= \varepsilon^k \cdot S_{k+1} n_{k+1}^c + O(\varepsilon^{k+1})} + \overbrace{L(\varepsilon) \cdot p_k(\varepsilon) \cdot n_{k+1}}^{= O(\varepsilon^{k+1})} \} \\ &\quad + \varepsilon^{2k-2} \cdot \mathcal{B}(\varepsilon) \cdot [\, p_k(\varepsilon) \cdot \{ \varepsilon^k \cdot n_1^c + \cdots + \varepsilon \cdot n_k^c + n_{k+1}^c + n_{k+1} \} ]^2 \\ &\quad + \varepsilon^{3k-3} \cdot r\left(\varepsilon, \varepsilon^{k-1} \cdot b\right) \cdot b^3 \end{aligned}$$

and as previously, we can simplify and collect powers of $\varepsilon$ according to

$$\begin{aligned} G[*] &= G[\, z(\varepsilon) \,] + \varepsilon^{2k-1} \cdot \{ S_1 \cdot n_1^c + \cdots + S_{k+1} \cdot n_{k+1}^c + O(\varepsilon) \} \\ &+ \varepsilon^{2k-2} \cdot [\, \mathcal{B}_0 + \varepsilon \cdot \mathcal{B}_1 + O(\varepsilon^2) \,] \cdot [\, (n_{k+1}^c + n_{k+1}) + \varepsilon \cdot (I_B \;\; \phi_1) \cdot \begin{pmatrix} n_k^c \\ n_{k+1}^c + n_{k+1} \end{pmatrix} + O(\varepsilon^2) \,]^2 \\ &+ \varepsilon^{3k-3} \cdot r\left(\varepsilon, \varepsilon^{k-1} \cdot b\right) \cdot b^3 \end{aligned} \tag{6.3}$$

$$\begin{aligned} &= G[\, z(\varepsilon) \,] + \varepsilon^{2k-1} \cdot \{ (S_1 \cdots S_{k+1}) \cdot n^c + O(\varepsilon) \} \\ &\qquad + \varepsilon^{2k-2} \cdot \mathcal{B}_0 \cdot (n_{k+1}^c + n_{k+1})^2 \\ &+ \varepsilon^{2k-1} \cdot [\, 2 \cdot \mathcal{B}_0 \cdot (n_{k+1}^c + n_{k+1}) \cdot (I_B \;\; \phi_1) \cdot \begin{pmatrix} n_k^c \\ n_{k+1}^c + n_{k+1} \end{pmatrix} + \mathcal{B}_1 \cdot (n_{k+1}^c + n_{k+1})^2 + O(\varepsilon) \,] \\ &\qquad + \varepsilon^{3k-3} \cdot r\left(\varepsilon, \varepsilon^{k-1} \cdot b\right) \cdot b^3 \end{aligned} \tag{6.4}$$

with $\mathcal{B}_1 = \frac{1}{2} G^{(3)}[0] \cdot z_1 \in L^2[B, \bar{B}]$. And again, our aim is to use the regularity of $(S_1 \cdots S_{k+1}) \in GL[N^c, \bar{B}]$ for linearization, requiring the coefficient to lowest power $\varepsilon^{2k-2}$ in (6.4) to be zero by

$$\mathcal{B}_0 \cdot (n_{k+1}^c + n_{k+1})^2 = 0\,, \tag{6.5}$$

which means zero curvature at $z = 0$ in the direction of the subspace $N_{k+1}^c \oplus N_{k+1}$. Then, the final form arises according to

$$\begin{aligned} G[*] = G[\,z(\varepsilon)\,] \;+\; \varepsilon^{2k-1} \cdot \{\,(S_1 \cdots S_{k+1}) \cdot n^c \\ +\; 2 \cdot \mathcal{B}_0 \cdot (n_{k+1}^c + n_{k+1}) \cdot (I_B \;\; \phi_1) \cdot \begin{pmatrix} n_k^c \\ n_{k+1}^c + n_{k+1} \end{pmatrix} \;+\; \mathcal{B}_1 \cdot (n_{k+1}^c + n_{k+1})^2\,\} \\ +\; \varepsilon^{2k} \cdot R(\,\varepsilon, n^c, n_{k+1}\,) \end{aligned} \tag{6.6}$$

with $R(\varepsilon, 0,0) = 0$. Note also, now we have to restrict to $k \geq 3$ for satisfying $3k - 3 \geq 2k$. Formula (6.6) corresponds to (3.5) in case of Jordan chains starting at $\varepsilon^k$, whereas the remainder equation (3.6) turns into

$$\begin{aligned} H(\varepsilon, n^c, n_{k+1}, \varphi) \coloneqq (S_1 \cdots S_{k+1}) \cdot (n^c - \varphi) \\ + 2 \cdot \mathcal{B}_0 \cdot (n_{k+1}^c + n_{k+1}) \cdot (I_B \;\; \phi_1) \cdot \begin{pmatrix} n_k^c \\ n_{k+1}^c + n_{k+1} \end{pmatrix} \;+\; \mathcal{B}_1 \cdot (n_{k+1}^c + n_{k+1})^2 \\ +\; \varepsilon \cdot R(\varepsilon, n^c, n_{k+1}) = 0\,. \end{aligned} \tag{6.7}$$

Again, we have $H(\varepsilon, 0,0,0) = 0$ and $H_{n^c}(0,0,0,0) = (S_1 \cdots S_{k+1}) \in GL[N^c, \bar{B}]$, yielding the existence of a locally unique and smooth operator function $\psi^c(\varepsilon, \varphi, n_{k+1})$ satisfying

$$H[\,\varepsilon,\;\; \psi^c(\varepsilon, \varphi, n_{k+1}),\;\; n_{k+1},\;\; \varphi\,] = 0 \tag{6.8}$$

and linearization is obtained according to (6.2), (6.6), (6.7) by

$$\begin{aligned} G[\,z(\varepsilon) \;+\; \varepsilon^{k-1} \cdot p_k(\varepsilon) \cdot \{\,\varepsilon^k \cdot \psi_1^c(\varepsilon, \varphi, n_{k+1}) \;+ \cdots +\; \psi_{k+1}^c(\varepsilon, \varphi, n_{k+1}) \;+\; n_{k+1}\,\}\,] \\ = G[\,z(\varepsilon)\,] \;+\; \varepsilon^{2k-1} \cdot (S_1 \cdots S_{k+1}) \cdot \varphi\,. \end{aligned} \tag{6.9}$$

Summarizing, for $k \geq 3$ and under the additional assumption (6.5), the linearization Theorem 1 is also valid for Jordan chains starting at $\varepsilon^{k-1}$, with red marked $\varepsilon^k$ and $\varepsilon^{2k}$ in Theorem 1 to be replaced by $\varepsilon^{k-1}$ and $\varepsilon^{2k-1}$ respectively. In particular, the size of the cone valid for linearization is increased by the order of 1.

In the next step, let us turn to the zeros of Corollary 2, but now with starting position of Jordan chains at $\varepsilon^{k-1}$ and assumed validity of the curvature condition (6.5). Then, by considering (6.6), we merely have to require $z(\varepsilon)$ to be an approximation of order $\varepsilon^{2k-1}$ according to

$$G[\,z(\varepsilon)\,] = \varepsilon^{2k-1} \cdot \bar{b}(\varepsilon) \tag{6.10}$$

for arriving at the remainder equation

$$\begin{aligned} (S_1 \cdots S_{k+1}) \cdot n^c \;+\; 2 \cdot \mathcal{B}_0 \cdot (n_{k+1}^c + n_{k+1}) \cdot (I_B \;\; \phi_1) \cdot \begin{pmatrix} n_k^c \\ n_{k+1}^c + n_{k+1} \end{pmatrix} \\ +\; \mathcal{B}_1 \cdot (n_{k+1}^c + n_{k+1})^2 \;+\; \bar{b}(\varepsilon) \;+\; \varepsilon \cdot R(\,\varepsilon, n^c, n_{k+1}\,) \end{aligned} \tag{6.11}$$

$$=: (S_1 \cdots S_{k+1}) \cdot n^c \; + \; \mathcal{F}(\, n_k^c, n_{k+1}^c, n_{k+1}\,) \; + \; \bar{b}(\varepsilon) \; + \; \varepsilon \cdot R(\, \varepsilon, n^c, n_{k+1}\,) = 0\,,$$

which replaces previous equation (5.2). Further, splitting equation (6.11) by subspaces of $\bar{B} = R_1 \oplus \cdots \oplus R_{k+1}$, we obtain at $\varepsilon = 0$

$$\begin{array}{lllllllll} S_1 \cdot n_1^c & + & \mathcal{P}_1 \cdot \mathcal{F}(\, n_k^c, n_{k+1}^c, n_{k+1}\,) & + & \mathcal{P}_1 \cdot \bar{b}(0) & = & 0 & & \\ & \vdots & & \vdots & & & & & \\ S_{k-1} \cdot n_{k-1}^c & + & \mathcal{P}_{k-1} \cdot \mathcal{F}(\, n_k^c, n_{k+1}^c, n_{k+1}\,) & + & \mathcal{P}_{k-1} \cdot \bar{b}(0) & = & 0 & & (6.12) \\ S_k \cdot n_k^c & + & \mathcal{P}_k \cdot \mathcal{F}(\, n_k^c, n_{k+1}^c, n_{k+1}\,) & + & \mathcal{P}_k \cdot \bar{b}(0) & = & 0 & & \\ S_{k+1} \cdot n_{k+1}^c & + & \mathcal{P}_{k+1} \cdot \mathcal{F}(\, n_k^c, n_{k+1}^c, n_{k+1}\,) & + & \mathcal{P}_{k+1} \cdot \bar{b}(0) & = & 0 & & \end{array}$$

with quadratic last two equations concerning $n_k^c$ and $n_{k+1}^c$. Again, within (6.12) decoupling of the equations into a quadratic and a linear part takes place, which requires now smallness of $\bar{b}(0)$ with respect to both of the subspaces $N_k^c$ and $N_{k+1}^c$ according to

$$\|\, \mathcal{P}_k \cdot \bar{b}(0) \,\| \ll 1 \quad \textit{and} \quad \|\, \mathcal{P}_{k+1} \cdot \bar{b}(0) \,\| \ll 1 \tag{6.13}$$

for solvability of the last two equations by $\bar{n}_k^c(0, n_{k+1}) \sim 0$ and $\bar{n}_{k+1}^c(0, n_{k+1}) \sim 0$ for $\|n_{k+1}\| \ll 1$. Here we use $S_k \in GL[N_k^c, R_k]$ and $S_{k+1} \in GL[N_{k+1}^c, R_{k+1}]$, as well as the Banach Lemma for control of the perturbation.

Finally, first $k-1$ equations in (6.12) are solved at $\varepsilon = 0$ by linear inversion, whereas the continuation to $\varepsilon \neq 0$ follows, as before, from bounded bijectivity of the linearization in (6.11) at these solutions, thereby using smallness of $n_{k+1}$ and, as previously, $(S_1 \cdots S_{k+1}) \in GL[N^c, \bar{B}]$.

Concerning Newton's Lemma with starting position of Jordan chains at $\varepsilon^{k-1}$, let us summarize the results.

**Corollary 3:** Assume the family $L(\varepsilon)$ to be $k$-surjective with $k \geq 3$ and all subspaces closed in (1.10). Further, let $z(\varepsilon)$ be an approximation of $G[z] = 0$ satisfying

$$G[\, z(\varepsilon)\,] = \varepsilon^{2k-1} \cdot \bar{b}(\varepsilon) \tag{6.14}$$

$$\|\, (\mathcal{P}_k + \mathcal{P}_{k+1}) \cdot \bar{b}(0) \,\| \ll 1 \tag{6.15}$$

$$\mathcal{B}_0 \cdot (n_{k+1}^c + n_{k+1})^2 = 0\,. \tag{6.16}$$

Then, in the cone

$$z = z(\varepsilon) + \varepsilon^{k-1} \cdot p_k(\varepsilon) \cdot \{\, \varepsilon^k \cdot n_1^c \; + \cdots + \; \varepsilon \cdot n_k^c \; + \; n_{k+1}^c + n_{k+1}\,\}\,, \tag{6.17}$$

there exist smooth mappings $\bar{n}_i^c(\varepsilon, n_{k+1}), \varepsilon \in U, \|\, n_{k+1} \,\| \ll 1, i = 1, \cdots, k+1$, such that

$$z = z(\varepsilon) + \varepsilon^{k-i} \cdot p_k(\varepsilon) \cdot \{\, \varepsilon^k \cdot \bar{n}_1^c(\varepsilon, n_{k+1}) + \cdots + \varepsilon \cdot \bar{n}_k^c(\varepsilon, n_{k+1}) + \bar{n}_{k+1}^c(\varepsilon, n_{k+1}) + n_{k+1}\,\} \tag{6.18}$$

defines all solutions of $G[z] = 0$ in the cone.

Compared to Corollary 2, the impact of higher order terms on the action of $G[z]$ is further weakened by requiring zero curvature in the origin by (6.16) with respect to $N_{k+1}^c \oplus N_{k+1}$. But then, the cone is allowed to be enlarged by the order of 1 according to (6.17) and the approximation requirement for $z(\varepsilon)$ can be relaxed, also by the order of 1, according to (6.14), (6.15).

Along these lines, we may obtain some sort of sequence of linearization theorems and adapted versions of Newton Lemmas by consecutive shift of Jordan chains to the left, thereby imposing appropriate solvability conditions concerning the resulting equations. A first example is given by Corollary 3 using a one step shift to the left.

A second example is summarized in Corollary 4 below. Here, a sequence of Newton Lemmas is stated, which is characterized by rather strict assumptions concerning derivatives of $G$ in the origin, yet allowing starting positions of the Jordan chains from $\varepsilon^{k-1}$ down to $\varepsilon^1$.

**Corollary 4:** Assume the family $L(\varepsilon)$ to be $k$-surjective with $k \geq 2$ and all subspaces closed in (1.10). Further, assume the existence of $i \in \{1, \cdots, k-1\}$ with

$$G^{(2)}[\,0\,] = 0\,, \;\cdots\;, G^{(i+1)}[\,0\,] = 0 \tag{6.19}$$

and let $z(\varepsilon)$ be an approximation of $G[z] = 0$ satisfying

$$G[\,z(\varepsilon)\,] = \varepsilon^{2k-i} \cdot \bar{b}(\varepsilon) \tag{6.20}$$

$$\| \,\mathcal{P}_{k+1} \cdot \bar{b}(0)\, \| \ll 1\,. \tag{6.21}$$

Then, in the cone

$$z = z(\varepsilon) + \varepsilon^{k-i} \cdot p_k(\varepsilon) \cdot \{\, \varepsilon^k \cdot n_1^c \;+ \cdots + \varepsilon \cdot n_k^c \;+\; n_{k+1}^c + n_{k+1} \,\}\;, \tag{6.22}$$

there exist smooth mappings $\bar{n}_i^c(\varepsilon, n_{k+1}), \varepsilon \in U, \| \,n_{k+1}\, \| \ll 1, i = 1, \cdots, k+1$, such that

$$z = z(\varepsilon) + \varepsilon^{k-i} \cdot p_k(\varepsilon) \cdot \{\, \varepsilon^k \cdot \bar{n}_1^c(\varepsilon, n_{k+1}) + \cdots + \varepsilon \cdot \bar{n}_k^c(\varepsilon, n_{k+1}) + \bar{n}_{k+1}^c(\varepsilon, n_{k+1}) + n_{k+1} \,\} \tag{6.23}$$

define all solutions of $G[z] = 0$ in the cone.

In case of $i = 1$, assumption (6.19) requires the curvature of $G$ at $z = 0$ to be zero with respect to all directions, i.e. compared to (6.16) of Corollary 3, further restrictions occur, cancelling completely the $\mathcal{B}_0$ term in (6.11), and thus slightly changing the discussion of the remainder equation. On the other hand, the approximation requirement (6.21) can be weakened compared to (6.15).

For $i \in \{2, \cdots, k-1\}$, further derivatives of $G$ at $z = 0$ are assumed to be zero by (6.19), allowing to decrease required order of approximation down to $\varepsilon^{2k-i}$ and to enlarge the cone according to $\varepsilon^{k-i}$ in (6.22). Then, presupposed approximation $z(\varepsilon)$ and obtained solutions agree up to the order of $\varepsilon^{k-i-1}$. The simple proof of Corollary 4 is omitted.

We end this section with an example concerning calculation of $k$-surjectivity with corresponding cones, as well as application of the topological degree to show the existence of secondary solution curves. In addition, Corollary 4 is applied and the Milnor number $\mu$ is calculated by $k$ values of different solution branches.

**Example:** Given the real polynomial equation

$$G[\,x, y\,] = -xy^3 + x^5 + y^5 = 0\,, \tag{6.24}$$

composed of monomials of order 4 and 5. Neglecting high order monomials $x^5$ and $y^5$ respectively, the simplified equations

$$-xy^3 + y^5 = y^3 \cdot (-x + y^2) = 0 \qquad and \qquad -xy^3 + x^5 = x \cdot \overbrace{(-y^3 + x^4)}^{E_6} = 0 \tag{6.25}$$

arise with solution curves, resulting from brackets equal to zero, given by

$$x(y) = y^2 \qquad and \qquad y(x) = x^{\frac{4}{3}} . \tag{6.26}$$

Note that the bracket within the second equation in (6.25) represents an $ADE$-singularity [1] of type $E_6$ with nonsmooth solution curve $y(x)$ considering its first derivative at $x = 0$. The solution curves $x(y)$ and $y(x)$ in (6.26) of the simplified equations (6.25) are tangentially touching $y$- and $x$-axis respectively, and hence, it is reasonable to start with $y$-axis and $x$-axis as first approximations for calculation of true solutions of the full equation (6.24).

First, the $y$-axis is considered, parametrized according to

$$z(\varepsilon) = \varepsilon \cdot \begin{pmatrix} 0 \\ 1 \end{pmatrix} , \tag{6.27}$$

which is now investigated with respect to $k$-surjectivity of the associated family

$$L(\varepsilon) = G'[\, z(\varepsilon)\, ] = [\, -y^3 + 5x^4 \mid -3xy^2 + 5y^4\, ]_{\substack{x=0\\ y=\varepsilon}} \tag{6.28}$$

$$= [\, -\varepsilon^3 \mid 0\, ] \in L[\, \mathbb{R}^2, \mathbb{R}\, ] .$$

Surjectivity arises with speed of order $k = 3$, yielding 3-surjectivity of $L(\varepsilon)$. We remark that in general the concrete calculation of the operators $S_1, \cdots, S_{k+1}$ simplifies extremly along a coordinate axis and in case of $G: B = \mathbb{R}^2 \longrightarrow \bar{B} = \mathbb{R}$. For $i \geq 1$, we obtain along the $y$-axis

$$S_i = c_i \cdot [\, G_{xy^{i-1}}[0] \quad G_{y^i}[0]\, ] \quad with \quad c_i \neq 0 , \tag{6.29}$$

i.e. in case of (6.24), we have $S_1 = S_2 = S_3 = [0 \;\; 0]$ and $S_4 = c_4 \cdot [G_{xy^3}[0] \;\; 0] = c_4 \cdot [-6 \;\; 0]$. Note that in case of $\bar{B} = \mathbb{R}$, surjectivity is obtained, as soon as a first entry different from zero occurs in $S_i \in L[\, \mathbb{R}^2, \mathbb{R}\, ]$.

Next, we have $G^{(2)}[0] = 0$ as well as $G^{(3)}[0] = 0$ and we can try to apply Corollary 4 with $k = 3$ and $i = 2$. Note also $N_{k+1} = N_4 = N[S_4] = span\, [0 \;\; 1]^T$ and $\mathcal{P}_{k+1} = \mathcal{P}_4 = I_{\bar{B}} = 1$.

Now, the $y$-axis obviously implies an approximation of order $\varepsilon^5$ and considering $2k - i = 4$ from (6.20), we end up with

$$G[\, z(\varepsilon)\, ] = G[\, 0, \varepsilon\, ] = \varepsilon^5 = \varepsilon^4 \cdot \varepsilon = \varepsilon^4 \cdot \bar{b}(\varepsilon) , \tag{6.30}$$

implying (6.21) with $\bar{b}(0) = 0$ and the existence of a nontrivial solution curve along the $y$-axis is shown by (6.23), (6.27) according to $z = z(\varepsilon) + O(\varepsilon^{k-i}) = [0 \;\; \varepsilon]^T + O(\varepsilon)$.

In figure 8 below, this solution curve is indicated in the left diagram by red colour. Other level sets (black) are also depicted in the cone of order $\varepsilon^{k-i} = \varepsilon^1$.

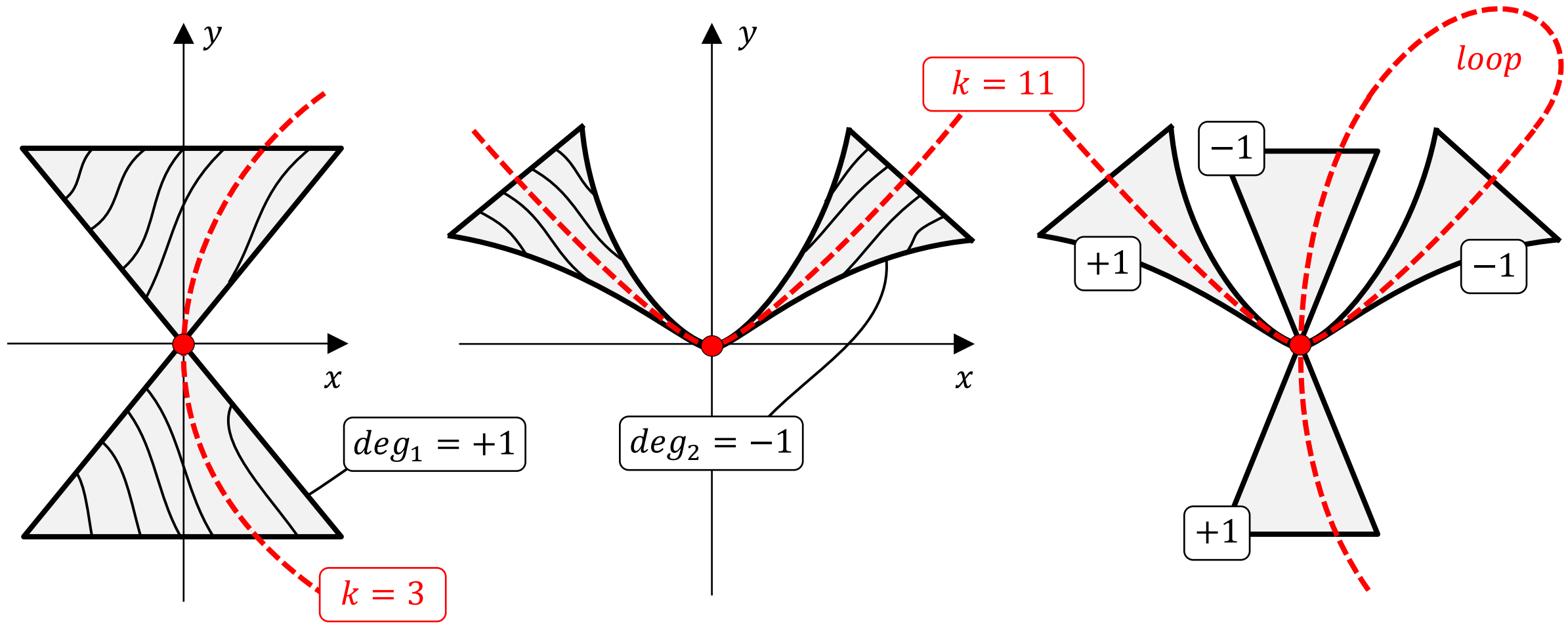

Figure 8: Solution curves (red) with corresponding cones (grey) of the singularity (6.24).

Note that standard Newton Lemma cannot be applied along the $y$-axis to ensure the existence of the solution curve, due to

$$det\ G_x[\,0,\varepsilon\,] = -\varepsilon^3$$

$$G[\,0,\varepsilon\,] = \varepsilon^5 \nleq (\,det\ G_x[\,0,\varepsilon\,]\,)^2 \cdot \varepsilon = \varepsilon^7\,. \tag{6.31}$$

Next, let us turn to the second solution curve, supposed to be positioned near the zero set of the $E_6$ singularity. First, performing a similar calculation as above, it is straightforward to see that even Corollary 4 cannot be applied to ensure the existence of this solution curve when starting from the $x$-axis and we are forced to improve the approximation by choosing the approximation $y(x)$ from (6.26) according to

$$y(x) = x^{\frac{4}{3}} \qquad \Leftrightarrow \qquad z(\varepsilon) = \begin{pmatrix} \varepsilon^3 \\ \varepsilon^4 \end{pmatrix}. \tag{6.32}$$

Also with this improvement, standard Newton lemma cannot be applied, because of

$$det\ G_y[\,\varepsilon^3,\varepsilon^4\,] \;=\; [\,-3xy^2 + 5y^4\,]_{\substack{x=\varepsilon^3 \\ y=\varepsilon^4}} \;=\; \varepsilon^{11}\cdot(\,-3+5\varepsilon^5\,)$$

$$G[\,\varepsilon^3,\varepsilon^4\,] \;=\; [\,-xy^3 + x^5 + y^5\,]_{\substack{x=\varepsilon^3 \\ y=\varepsilon^4}} \;=\; [\,y^5\,]_{y=\varepsilon^4} = \varepsilon^{20}\,, \tag{6.33}$$

implying

$$G[\,\varepsilon^3,\varepsilon^4\,] = \varepsilon^{20} \nleq \left(det\ G_y[\,\varepsilon^3,\varepsilon^4\,]\right)^2 \cdot \varepsilon = \;\varepsilon^{23}\cdot(\,-3+5\varepsilon^5\,)^2\,. \tag{6.34}$$

Now, concerning the application of Corollary 4 with $z(\varepsilon)$ from (6.32), we obtain

$$L(\varepsilon) = G'[\,\varepsilon^3,\varepsilon^4\,] = [\,-y^3 + 5x^4 \;\mid\; -3xy^2 + 5y^4\,]_{\substack{x=\varepsilon^3 \\ y=\varepsilon^4}} \tag{6.35}$$

$$= [\,\varepsilon^{12}\cdot 4 \;\mid \varepsilon^{11}\cdot(\,-3+5\varepsilon^5\,)\,]$$

and surjectivity of $L(\varepsilon)$ occurs with speed of order $k = 11$, yielding 11-surjectivity of $L(\varepsilon)$. Concrete calculation of the operators $S_i$ implies $S_1 = \cdots = S_{11} = 0$ as well as $S_{12} \neq 0$ and condition (6.20) in Corollary 4 is satisfied with $k = 11$ and $i = 2$ by

$$G[\,z(\varepsilon)\,] = G[\,\varepsilon^3, \varepsilon^4\,] = \varepsilon^{20} \cdot 1 = \varepsilon^{2\cdot 11-2} \cdot 1 = \varepsilon^{20} \cdot \bar{b}(\varepsilon)\,. \tag{6.36}$$

However, due to $\bar{b}(\varepsilon) = 1$, condition (6.21) is not yet satisfied, unless we force smallness of $\bar{b}(\varepsilon)$ by a parameter $\alpha$ according to

$$G[\,x, y\,] = -xy^3 + x^5 + \alpha \cdot y^5 = 0 \tag{6.37}$$

with $|\alpha| \ll 1$. Then, (6.33) turns into

$$G[\,\varepsilon^3, \varepsilon^4\,] = [\,-xy^3 + x^5 + \alpha \cdot y^5\,]_{\substack{x=\varepsilon^3 \\ y=\varepsilon^4}} = [\,\alpha \cdot y^5\,]_{y=\varepsilon^4} = \varepsilon^{20} \cdot \alpha \tag{6.38}$$

and (6.21) is assured too, yielding by (6.23) with $\varepsilon^{k-i} = \varepsilon^9$ and (6.32) the existence of a smooth solution curve of the form

$$z(\varepsilon) \;=\; \varepsilon^3 \cdot \begin{pmatrix} 1 \\ 0 \end{pmatrix} \;+\; \varepsilon^4 \cdot \begin{pmatrix} 0 \\ 1 \end{pmatrix} \;+\; O(\,\varepsilon^9\,)\,, \tag{6.39}$$

locally parametrized by the external parameter $\varepsilon \in U$. The corresponding constellation is depicted in the middle diagram of figure 8 and all in all, the two solution curves along $y$- and $x$-axis are globally connected by a loop, as shown in the diagram on the right hand side.

It is also interesting to note that each of the two solution curves can be taken as a first curve for proving the existence of the other solution curve by means of the topological degree. For this purpose, we remark that each of the half cones in figure 8 possesses a constant topological degree, which can be calculated by

$$deg_1(\varepsilon) \coloneqq sign\{\, det\, G_x[\,0, \varepsilon\,]\,\} = sign\{\, -\varepsilon^3\,\} = \pm 1 \tag{6.40}$$

$$deg_2(\varepsilon) \coloneqq sign\big\{\, det\, G_y[\,\varepsilon^3, \varepsilon^4\,]\,\big\} \;=\; sign\{\, \varepsilon^{11} \cdot (\,-3 + 5\varepsilon^5\,)\,\} = \pm 1\,. \tag{6.41}$$

Hence, different signs of the topological degree occur in the two half cones associated to each of the solution curves, implying a continuum of secondary solutions to emanate from the origin, this continuum producing the other solution curve. For details and applications in higher dimensions, see [17].

Finally, note that the Milnor number $\mu$ of the singularity (6.24) is easily calculated by the $k$ values $k_1 = 3$ and $k_2 = 11$ of the two solution curves through the origin according to

$$\mu = k_1 + k_2 - ord(G) + 1 = 3 + 11 - 4 + 1 = 11\,, \tag{6.42}$$

which represents an application of the general formula

$$\mu = k_1 + \cdots + k_\tau - ord(G) + 1 \tag{6.43}$$

for calculation of the Milnor number $\mu$, as derived in [16]. The formula is valid, if all segments of the Newton polygon associated to a singularity factorize completely with multiplicity 1. Here, $\tau$ denotes the number of different solution curves through the singularity. See also [16] for general application of (6.43) to simple $ADE$-singularities and for a possible generalization of (6.43) with respect to higher multiplicities of the segments.

Concerning the application of the topological degree, we note that formulas (6.40), (6.41) are only valid, if each crossing section contains exactly one solution element. This is one of the reasons, why it seems appropriate to show in the last section of this paper, how to eliminate the redundancy present within the parametrizations of Theorem 1 and subsequent Corollaries. In particular, we derive unique representations of the solutions of $G[z] = 0$.

## *7. Elimination of Redundancy*

For simplicity, this final section restricts to the starting position (3.1), (3.2) of the Jordan chains at $\varepsilon^k$, as used in Theorem 1 and Corollaries 1, 2.

Then by construction, the $\varepsilon$-parametric family of local diffeomorphisms $z = z(\varepsilon) + \varepsilon^k \cdot p_k(\varepsilon) \cdot \{\Psi^c(\varepsilon, \varphi, n_{k+1}) + n_{k+1}\}$ in (1.17) shows a strong redundancy by the fact that every element $b$ in the cone can be represented by each diffeomorphism of the family with image containing $b$. In some more detail, for given $b$ in the cone, there exist, by implicit function theorem, smooth functions $\varphi(b, \varepsilon) \in N^c$ and $n_{k+1}(b, \varepsilon) \in N_{k+1}$ satisfying

$$b = z(\varepsilon) + \varepsilon^k \cdot p_k(\varepsilon) \cdot \{\, \Psi^c[\, \varepsilon, \varphi(b, \varepsilon), n_{k+1}(b, \varepsilon)\,] + n_{k+1}(b, \varepsilon)\,\}\,. \tag{7.1}$$

Compare also the diagram in figure 1 on the right hand side.

In case of $k = 0$, this redundancy may simply be omitted by choosing only one of the diffeomorphisms from (1.17), e.g. the diffeomorphism that arises for $\varepsilon = 0$

$$z = z(0) + \varepsilon^0 \cdot p_k(0) \cdot \{\, \Psi^c(0, \varphi_1, n_1) + n_1\,\} = \Psi^c(0, \varphi_1, n_1) + n_1\,. \tag{7.2}$$

Now, in case of $k \geq 1$, the $\varepsilon$-redundancy in (7.1) can be removed by imposing an appropriate parametrization condition guided by the leading order term $\varepsilon^l \cdot z_l$ of the center curve

$$z(\varepsilon) \;=\; \varepsilon^l \cdot \underbrace{z_l}_{\neq 0} \;+ \cdots +\; \varepsilon^k \cdot z_k \;+ \cdots\,, \tag{7.3}$$

satisfying $1 \leq l \leq k$ for $k$-surjectivity to hold true. Then for every $b$ in the cone, exactly one of the possible $\varepsilon$ values in (7.1) may be sorted out by the following case-by-case analysis.

If $z_l \notin N^c$, then the projection $z_{l,k+1} := P_{N_{k+1}} z_l$ of $z_l$ to $N_{k+1}$ is different from 0 and the direct sum of $B$ in (1.10) can be refined by

$$B \;=\; N^c \,\oplus\, N_{k+1} \;=\; N^c \,\oplus\, \Pi_{k+1} \,\oplus\, \{\, z_{l,k+1}\,\}\,, \tag{7.4}$$

and again we require the complement $\Pi_{k+1}$ to be closed. Now, the restriction of $N_{k+1}$ to $\Pi_{k+1}$ within the family of diffeomorphisms (1.17) reads

$$z \;=\; z(\varepsilon) + \varepsilon^k \cdot p_k(\varepsilon) \cdot \{\, \Psi^c(\, \varepsilon, \varphi, p_{k+1}\,) + p_{k+1}\,\} \tag{7.5}$$

with $p_{k+1} \in \Pi_{k+1}$ and center curve $z(\varepsilon)$ given by

$$z(\varepsilon) \;=\; \varepsilon^l \cdot z_l \;+ \cdots =\; \varepsilon^l \cdot (\, \underbrace{z_l^c}_{\in N^c} \;+\; \underbrace{z_{l,k+1}}_{\neq 0}\,) \;+ \cdots\,. \tag{7.6}$$

Further, for fixed $\varepsilon \neq 0$, the mapping $p_k(\varepsilon) \cdot \{\Psi^c(\varepsilon, \varphi, p_{k+1}) + p_{k+1}\}$ is a local diffeomorphism within the subspace $N^c \oplus \Pi_{k+1}$, which represents by (7.4) a direct complement to the subspace spanned by $z_{l,k+1}$. Hence, at every point of the center line $z(\varepsilon)$, the set $\varepsilon^k \cdot p_k(\varepsilon) \cdot \{\Psi^c(\varepsilon, \varphi, p_{k+1}) + p_{k+1}\}$ from the complement $N^c \oplus \Pi_{k+1}$ is attached, implying (7.5) to define a mapping from the domain

$$(\varepsilon, \ \varphi, \ p_{k+1}) \in$$

$$\hat{V} := U \setminus \{0\} \times \{ \varphi \in N^c \mid \|\varphi\| \ll 1 \} \times \{ p_{k+1} \in \Pi_{k+1} \mid \|p_{k+1}\| \ll 1 \} \qquad (7.7)$$

to an open cone around $z(\varepsilon)$ in $B$ with the property that for given $(\bar{\varepsilon}, \bar{\varphi}, \bar{p}_{k+1}), \bar{\varepsilon} \neq 0$ mapping (7.5) represents a local diffeomorphism between an open neighbourhood of $(\bar{\varepsilon}, \bar{\varphi}, \bar{p}_{k+1})$ and a corresponding neighbourhood of $z(\bar{\varepsilon}) + \bar{\varepsilon}^k \cdot p_k(\bar{\varepsilon}) \cdot \{\psi^c(\bar{\varepsilon}, \bar{\varphi}, \bar{p}_{k+1}) + \bar{p}_{k+1}\}$ in $B$.

In some more detail, the proof resembles the proof of Theorem 1 (ii), with the difference that the derivative in (3.24) with respect to $N_{k+1}$ is replaced by derivatives with respect to $\Pi_{k+1} \subsetneq N_{k+1}$ and $\varepsilon$. This is possible, due to the assumption $z_l \notin N^c$, implying

$$z'(\varepsilon) \ = \ l \cdot \varepsilon^{l-1} \cdot \left( z_l^c + z_{l,k+1} \right) + \cdots \ \neq 0 \qquad (7.8)$$

for $\varepsilon \neq 0$. It remains to look at the case of $z_l \in N^c$. First, if $N_{k+1} \neq \{0\}$, then it is not difficult to see that the complement $N^c$ can slightly be perturbed for ensuring $z_l \notin N^c$ and we can return to (7.4). For details see [17].

However, if $N_{k+1} = \{0\}$, we have $B = N^c$ and $z_l \in N^c$ cannot be omitted. Now, as a consequence of $N_{k+1} = \{0\}$, a straightforward calculation shows that the leading order $l$ of $z(\varepsilon)$ will agree with the order $k$ of surjectivity and the $\varepsilon$-parametric family of diffeomorphisms (1.17) turns into

$$z \ = \ z(\varepsilon) + \varepsilon^k \cdot p_k(\varepsilon) \cdot \Psi^c(\varepsilon, \varphi) \qquad (7.9)$$

$$\stackrel{l=k}{=} \ \underbrace{\varepsilon^k \cdot (z_k + \cdots)}_{curve} \ + \ \underbrace{\varepsilon^k \cdot p_k(\varepsilon) \cdot \Psi^c(\varepsilon, \varphi)}_{image\ blow\ up}$$

with the effect that the speed of movement $\varepsilon^k$ along the curve $z(\varepsilon)$ equals the blow up speed $\varepsilon^k$ of the images $p_k(\varepsilon) \cdot \Psi^c(\varepsilon, \varphi)$. Then, we typically obtain a constellation, as depicted in the right diagram of figure 9 below.

For each $\varepsilon \neq 0$, the image of the diffeomorphism contains the origin (black dot) and we can fix an arbitrary $\varepsilon \neq 0$ in (7.9) for describing the behaviour of $G[z]$ in an open neighbourhood of $z = 0$, i.e. in this case also a single diffeomorphism can be used, comparable to the case of $k = 0$ from (7.2). In the right diagram of figure 9, the chosen diffeomorphism is indicated by grey colour.

All in all, working with the $\varepsilon$-parametric family of diffeomorphisms $z = z(\varepsilon) + \varepsilon^k \cdot p_k(\varepsilon) \cdot \{\Psi^c(\varepsilon, \varphi, n_{k+1}) + n_{k+1}\}$ from (1.17) offers the possibility to avoid the different cases from above.

On the other hand, describing the action of $G[z]$ without redundancy of charts, will allow us to create unique representations of the level sets in the cone and in particular, will supply us with a unique representation of the solutions of $G[z] = 0$.

In some more detail, for $k = 0$ the diffeomorphism from (7.2) and Theorem 1 (i) imply

$$G[\,\Psi^c(0, \varphi, n_1) + n_1\,] = G[\,0\,] \ + \ G'[0] \cdot \varphi \qquad (7.10)$$

and from

$$G[\,0\,] \;+\; G'[0]\cdot\varphi = \bar{b} \tag{7.11}$$

$$\Leftrightarrow \quad \varphi(\,\bar{b}\,) \coloneqq G'[0]^{-1}\cdot(\,\bar{b} - G[\,0\,]\,)\,,$$

we obtain a unique level set to the level $\bar{b} \in \bar{B}$ by

$$b(\,\bar{b}, n_1\,) \coloneqq \psi^c(\,0,\;\; \varphi(\bar{b}),\;\; n_1\,) \;+\; n_1\,. \tag{7.12}$$

In particular, the level set to the value $\bar{b} = G[0]$ reads $b(G[0], n_1) = \Psi^c(0,0,n_1) + n_1$. This configuration is indicated in the left diagram of figure 9.

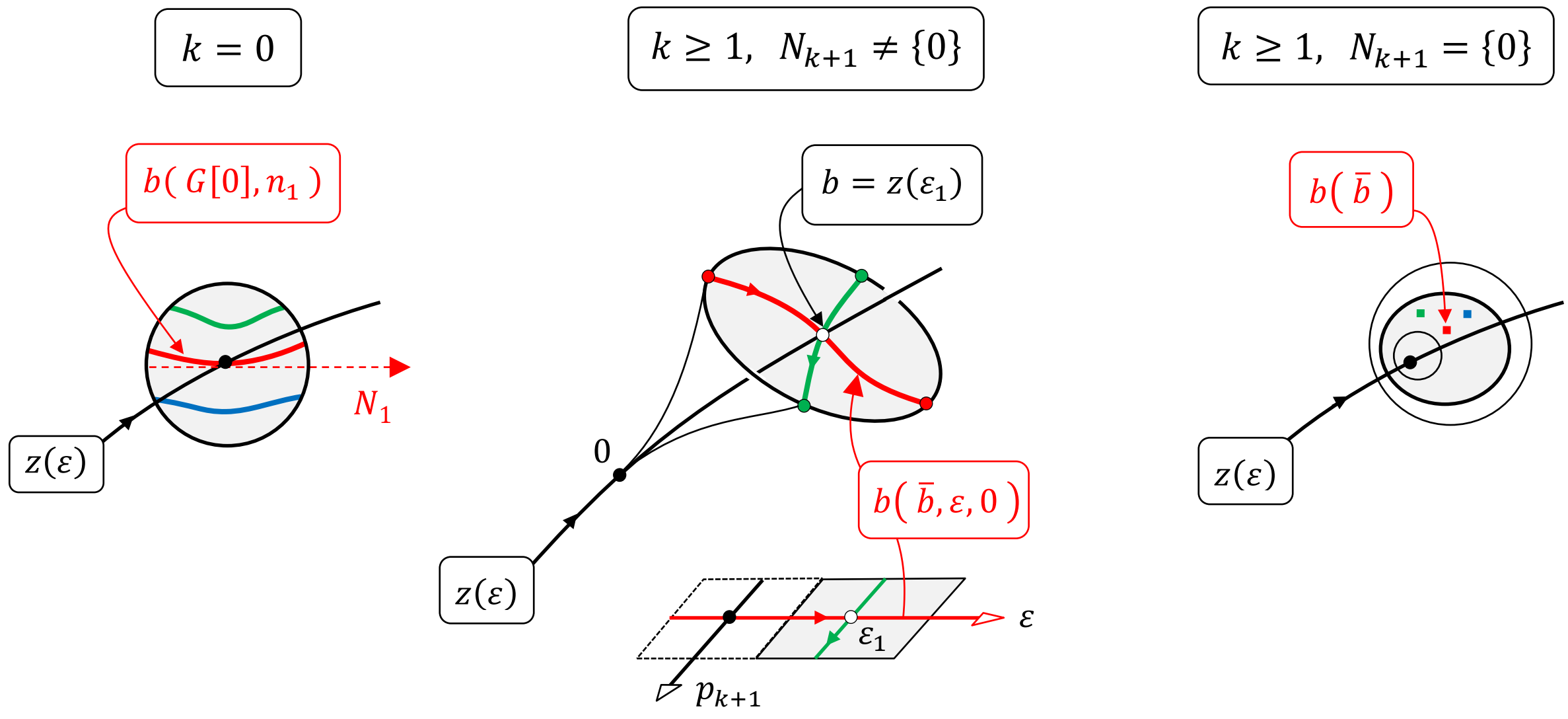

Figure 9: Different cases concerning level sets and elimination of redundancy.

Concerning the middle diagram with $k \geq 1$ and $N_{k+1} \neq \{0\}$, the local diffeomorphisms from (7.5) can be used to obtain from Theorem 1 (i)

$$G[\,z(\varepsilon) + \varepsilon^k\cdot p_k(\varepsilon)\cdot\{\,\Psi^c(\,\varepsilon, \varphi, p_{k+1}\,) + p_{k+1}\,\}\,]$$

$$= G[\,z(\varepsilon)\,] \;+\; \varepsilon^{2k}\cdot(\,S_1\cdots S_{k+1}\,)\cdot\varphi \tag{7.13}$$

and from

$$G[\,z(\varepsilon)\,] \;+\; \varepsilon^{2k}\cdot(\,S_1\cdots S_{k+1}\,)\cdot\varphi = \bar{b}$$

$$\Leftrightarrow \quad \varphi(\,\bar{b}, \varepsilon\,) \coloneqq \varepsilon^{-2k}\cdot(\,S_1\cdots S_{k+1}\,)^{-1}\cdot(\,\bar{b} - G[\,z(\varepsilon)\,]\,)\,, \tag{7.14}$$

we see that the level set to the value $\bar{b} \in \bar{B}$ is uniquely obtained by

$$b(\,\bar{b}, \varepsilon, p_{k+1}\,) \coloneqq z(\varepsilon) + \varepsilon^k\cdot p_k(\varepsilon)\cdot\{\,\Psi^c(\,\varepsilon,\;\; \varphi(\bar{b},\varepsilon),\;\; p_{k+1}\,) + p_{k+1}\,\}\,. \tag{7.15}$$

Due to $\varepsilon^{-2k}$ in (7.14) and $G[0] = 0$ for $k \geq 1$, the level set to a value $\bar{b} \neq 0$ cannot be extended up to $\varepsilon = 0$, but it is defined for all $(\varepsilon, p_{k+1})$ satisfying $\left(\varphi(\bar{b},\varepsilon), p_{k+1}\right) \in \hat{V}$ from (7.7). By (7.15), the level set is given in a maximal way up to the boundary of $\hat{V}$.

In the middle diagram of figure 9, the level set through the point $b = z(\varepsilon_1)$ is indicated by grey colour with corresponding value given by $\bar{b} = G[z(\varepsilon_1)]$. To show the effect of the level set parametrization by $(\varepsilon, p_{k+1})$, the Banach space $B$ is represented as $B = \mathbb{R}^3$ and the level sets turn into 2-dimensional isosurfaces when assuming $\bar{B} = \mathbb{R}$.

The red line shows the $\varepsilon$-coordinate line for fixed $p_{k+1} = 0$ within the isosurface through $z(\varepsilon_1)$. The green line indicates the $p_{k+1}$-coordinate line for fixed $\varepsilon = \varepsilon_1$. Finally, in the lower part of the middle diagram, only the grey shaded region in $(\varepsilon, p_{k+1})$-space can be used for parametrization of the isosurface, due to blow up of $\varphi(\bar{b}, \varepsilon)$ for $\bar{b} \neq 0$ and $\varepsilon \to 0$ in (7.14).

On the other hand, if we choose $\bar{b} = 0$ in (7.14) and require an approximation of order $2k$ according to $G[z(\varepsilon)] = \varepsilon^{2k} \cdot \bar{b}(\varepsilon)$, we obtain

$$\begin{aligned} \varphi(0, \varepsilon) &= -\varepsilon^{-2k} \cdot (S_1 \cdots S_{k+1})^{-1} \cdot G[z(\varepsilon)] \\ &= -\varepsilon^{-2k} \cdot (S_1 \cdots S_{k+1})^{-1} \cdot \varepsilon^{2k} \cdot \bar{b}(\varepsilon) \\ &= -(S_1 \cdots S_{k+1})^{-1} \cdot \bar{b}(\varepsilon) \end{aligned} \tag{7.16}$$

and the solutions of $G[z] = 0$ in the cone are uniquely given by

$$b(0, \varepsilon, p_{k+1}) = z(\varepsilon) + \varepsilon^k \cdot p_k(\varepsilon) \cdot \{ \Psi^c(\varepsilon, \varphi(0, \varepsilon), p_{k+1}) + p_{k+1} \}, \tag{7.17}$$

which are well defined for $\varepsilon \in U$ and $\|p_{k+1}\| \ll 1$. In addition, $\varphi(0, \varepsilon)$ has to be small, which is guaranteed by $\|\bar{b}(0)\| \ll 1$ in (1.24). In the punctured neighbourhood $\varepsilon \in U \backslash \{0\}$, the zeros $b(0, \varepsilon, p_{k+1})$ define a regular Banach manifold in $B$, as shown in [16].

Finally, in case of $k \geq 1$ and $N_{k+1} = \{0\}$, formula (7.15) turns into

$$b(\bar{b}) \coloneqq z(\varepsilon) + \varepsilon^k \cdot p_k(\varepsilon) \cdot \Psi^c(\varepsilon, \varphi(\bar{b}, \varepsilon)) \tag{7.18}$$

and by (7.9), the mapping $b(\bar{b})$ represents a diffeomorphism between neighbourhoods of $0 \in \bar{B}$ and $0 \in B$ for $\varepsilon \neq 0$, implying point level sets in $B$, as indicated in the right diagram of figure 9 by coloured dots.

**Acknowledgements:** The main ideas of the paper go back to E. Bohl, W.-J. Beyn and E. Jäger. In addition, the paper was motivated by stimulating discussions with J. López-Gómez. The author wishes to express his thanks to those people.

**Declarations**

**Ethical Approval:** not applicable

**Funding:** not applicable

**Availability of data and materials:** not applicable

References

[1] V. I. Arnold, S. M. Gusein-Zade, A. N. Varchenko, *Singularities of Differentiable Maps,* Volume **I**, Birkhäuser-Verlag, 1985.

[2] H. Bart, M.A. Kaashoek, D.C. Lay, *Stability properties of finite meromorphic operator functions I, II, III*, Indagationes Mathematicae, Volume **77**, Issue 3, 217-259, 1974.

[3] R. F. Coleman, H. J. Stein, *On Newton's lemma,* Journal of Algebra 322, 3427-3450, 2009.

[4] C. Bruschek, H. Hauser, *Arcs, Cords and Felts – Six Instances of the Linearization Principle,* arXiv:1006.5295v1 [math.AG], 2010.

[5] J. Esquinas, *Optimal multiplicity in local bifurcation theory, II: General case,* J. Diff. Equ. **75**, 206-215, 1988.

[6] B. Fisher, *A note on Hensel's Lemma in several variables,* Proc. Amer. Math. Soc. **125**, no. 11, 3185-3189, 1997.

[7] I. Gohberg, J. Leiterer, *Holomorphic Operator Functions of One Variable and Applications,* Operator Theory **192**, ISBN 978-3-0346-0125-2, Birkhäuser, 2009.

[8] M. J. Greenberg, *Rational points in Henselian discrete valuation rings,* Publ. Math. Inst. Hautes Études Sci. **31**, 59-64, 1966.

[9] H. Hauser, *The classical Artin Approximation theorems*, Bull. Amer. Math. Soc. **54**, 595-633, 2017.

[10] W. Kaballo, *Meromorphic generalized inverses of operator functions,* Indagationes Mathematicae **23**, 970-994, Elsevier, 2012.

[11] J. López-Gómez, C. Mora-Corral, *Algebraic Multiplicity of Eigenvalues of Linear Operators,* Operator Theory **177**, ISBN 978-3-7643-8400-5, Birkhäuser, 2007.

[12] J. López-Gómez, *Spectral Theory and Nonlinear Functional Analysis,* Research Notes in Mathematics **426**, Chapman and Hall/CRC Press, Boca Raton, 2001.

[13] R. Magnus, C. Mora-Corral, *Natural Representations of the Multiplicity of an Analytic Operator-valued Function at an Isolated Point of the Spectrum*, Integr. equ. oper. theory **53** (2005), 87–106, Birkhäuser, 2005.

[14] G. Rond, *Artin Approximation*, arXiv:1506.04717v4 [math.AC], 2018.

[15] M. Stiefenhofer, *Singular Perturbation and Bifurcation in case of Dictyostelium discoideum,* Ph.D. thesis, University of Konstanz (Germany), Research Gate M. Stiefenhofer, Hartung-Gorre-Verlag Konstanz, 1995.

[16] M. Stiefenhofer, *Direct sum condition and Artin Approximation in Banach spaces*, arXiv: 1905.07583v1 [math.AG], 2019.

[17] M. Stiefenhofer, *Fine Resolution of k-transversal Cones,* arXiv:2004.08615, [math. AG], 2020.

[18] M. Stiefenhofer, *Formal and Analytic Diagonalization of Operator functions,* arXiv:2305. 14228v1 [math.AG], 2023, submitted.

[19] M. Stiefenhofer, *A note on an expansion formula with application to nonlinear differential algebraic equations,* arXiv:1912.10735 [math.AG], 2019.